\documentclass[11pt]{article}
\usepackage{graphics}
\usepackage{latexsym}
\usepackage{amsfonts}

\usepackage{caption}
\usepackage{subcaption}
\usepackage[demo]{graphicx}
\usepackage{graphicx}
\usepackage{amsmath}
\usepackage{amsthm}
\usepackage{amssymb}
\newtheorem{definition}{Definition}
\newtheorem{theorem}{Theorem}

\newtheorem{assumption}{Assumption}

\newtheorem{remark}{Remark}

\newcommand*{\QEDB}{\hfill\ensuremath{\square}}

\begin{document}
	
	%First title page - only necessary on longer papers
	\begin{titlepage}
		\Large
		\begin{center}
			{\bf \Huge
				An Extremum-Seeking Co-Simulation Based Framework for Passivation Theory and its Application in Adaptive Cruise Control Systems
			}
		\end{center}
		
		\begin{center}
			Technical Report of the ISIS Group\\
			at the University of Notre Dame\\
			ISIS-2015-003\\ %fix report number
			August 2015 %fix date
		\end{center}
		\vspace{1.5in}
		
		\begin{center}
			Arash Rahnama$^1$, Meng Xia$^1$, Shige Wang$^2$ and Panos J. Antsaklis$^1$\\
			$^1$Department of Electrical Engineering\\
			University of Notre Dame\\
			Notre Dame, IN 46556\\
			$^2$General Motors R\&D,\\
			Warren, MI 48090\\
		\end{center}
		\vspace{1in}
		\begin{center}
			{\bf \Large Interdisciplinary Studies in Intelligent Systems}
		\end{center}
		\vspace{0.5in}
		\begin{center}
			{\textbf {Acknowledgements}    The support of the National Science Foundation under the CPS Grant No. CNS-1035655 is gratefully acknowledged.
			}
		\end{center}
	\end{titlepage}
	
%%	\title
%%	{An Extremum-Seeking Co-Simulation Based Framework for Passivation Theory and its Application in Adaptive Cruise Control Systems\vspace{0.2in}\newline
%%		\author{A. Rahnama, M. Xia, S. Wang and P. J. Antsaklis
			%%\thanks{Rahnama, Xia and Antsaklis are with the Department of Electrical Engineering, University of Notre
	%%			Dame, Notre Dame, IN 46556, USA (e-mail:~\{arahnama, ~mxia,~antsaklis.1\}~%~~mmccour1
	%%			@nd.edu). Wang is with General Motors R\&D (e-mail:~shige.wang@gm.com).
	%%			The support of the National Science Foundation under the CPS Grant No. CNS-1035655 is gratefully acknowledged. }}
%	}
%%	\maketitle
	%============================================================================================
	\begin{abstract}
		In this report, we apply an input-output transformation passivation method, described in our previous works, to an Adaptive Cruise Control system. We analyze the system's performance under a co-simulation framework that makes use of an online optimization method called extremum-seeking to achieve the optimized behavior. The matrix for passivation method encompasses commonly used methods of series, feedback and feed-forward interconnections for passivating the system. We have previously shown that passivity levels can be guaranteed for a system using our passivation method. In this work, an extremum-seeking algorithm was used to determine the passivation parameters.  It is known that systems with input-output time-delays are not passive. On the other hand, time-delays are unavoidable in automotive systems and commonly emerge in software implementations and communication units as well as driver's behavior. We show that by using our passivation method, we can passivate the system and improve its overall performance. Our simulation examples in CarSim and Simulink will show that the passive system has a considerably better performance.
	\end{abstract}
	\section{Introduction}
	\subsection{Cyber-Physical Systems}
	The necessity for designing smart compositional systems that are able to perform different control tasks in today technologically advanced world is increasing. Consequently, the interconnection amongst these units becomes important as well. Hence the accurate design of compositional systems that govern these interconnections is important. Each unit under this compositional design should perform as expected, and without the fear of any intervention that might disrupt the system's overall performance. This increased reliance on technological advancement accompanied with new advancements in sensing, communications, control and computation have created an emerging class of complex systems called Cyber-Physical Systems.
	
	Cyber-physical Systems (CPS) consist of large number of complex yet closely interconnected and integrated units. Examples of CPS may be found in smart transportation systems, smart medical devices, and smart energy systems \cite{antsaklis2012cyber}. Loosely speaking, such systems consist of two primary units: the physical part, which provides the system with a continuous model of the physical world, typically using ordinary differential equations, and a communication and computational part, which monitors, coordinates, and controls the physical systems. The computational unit includes the software component of the design and requires strong communication links in order to both receive and transfer the data to the physical world \cite{shang2013case}. The control and maintaining the robustness and reliability of cyber-physical systems presents huge challenges, The energy-based concepts of dissipativity and, more specifically, passivity provide a powerful tool for meeting the challenges that compositional systems produce. Passive systems do not generate energy, but only store or release the energy, which was provided to them. Under right conditions, passivity can present asymptotic stability for zero state detectable (ZSD) systems \cite{bao2007process}. Additionally, both negative feedback and parallel interconnections of passive systems stay passive, which means passivity and stability are preserved for large-scale systems consisting of passive stable units \cite{hill1976stability}. This renders passive designs a suitable candidate for designing cyber-physical systems.
	
	The underlying premise of our work is that an automobile meets all the criteria for a cyber-physical system. Automotive systems provide a platform consisting of both physical and computational components, both integrated through communication networks. Due to lack of a clear understanding of the complex and tight interactions that result from the integration of different control components, the design of automotive control applications is a complicated problem. Passivity has been studied as a possible systematic solution to this problem. This approach is beneficial in regards to improving overall system performance, and potentially can help solve the scalability challenges given the increasing integration of more and more control functionalities in vehicles. Due to the presence of delay, automotive systems are not passive. Our passivation method can passivate systems with input-output delays, hence our work introduces passivity as a useful tool for controlling sub-units in automotive systems.

	\label{sec:intr}

	\subsection{Systems with Input-Output Delay}
	%===================================================
	Computational delays, input delays, and measurement delays are some examples of delays in dynamical systems \cite{wang1998finite}. In addition approximation of higher order systems by lower order systems plus a time delay is a common approach in solving many control problems.  Delays usually introduce poor performance or instability to the closed loop system. In 1950’s, Otto Smith introduced a unique predictive control method that compensated for delay outputs using input values stored over a certain time window, and estimating the plant output accordingly \cite{smith1959controller}. Later this method combined with finite time integrals of the delayed input values was expanded to include unstable plants as well \cite{jang2009adaptive}.  Different adaptive control methods for handling time delays were developed later \cite{niculescu2003adaptive,kim2007application}.  Ortega and Lozano introduced a control version of delay systems, which is capable of handling the uncertainty and accumulated error that is usually the result of predictive control methodologies \cite{ortega1988globally}.  Many examples of systems including time delays in chemical, biological, mechanical, physiological, and electrical systems are given in \cite{kolmanovskii2013introduction,niculescu2001delay}, and a more detailed surveys of time delay systems can be found in \cite{richard2003time,gu2003survey}.

	\subsection{Adaptive Cruise Control}
	
	The Vanderbilt University model-based design of adaptive cruise control, explained in more details in \cite{eyisi2013model} was adapted as the main framework where our work was tested on. In our work the passivation method, and passivity control are applied to the development of an adaptive cruise control. An adaptive cruise control's main task is to maintain the desired velocity set by the driver while preserving a safe distance from other cars. The ACC accepts the speed of the lead vehicle, and the distance between the two cars as its inputs. The ACC will maintain the desire speed for the host car, as long as the safe distance between the host car and lead car is preserved, otherwise the velocity of the host car should be adjusted accordingly. This means that our system will have a hybrid structure of two modes, one to control vehicle when it is accelerating, and one for when it is decelerating. 
	
	The design of the adaptive cruise control (ACC) has been studied in the literature. Multiple-surface sliding control containing an “upper” sliding control, and a switching rule between brake or throttle control was used in \cite{gerdes1997vehicle, hedrick2000multiple} which results in a smooth alternation between both controllers by solely relying on the vehicle state and avoiding high frequency oscillation and risk of competing control inputs. Optimal control was used in \cite{ioannou1993autonomous, yi2001investigation} to design an ACC, and the performance was compared with that of the human driver models. The results are that the performance is much smoother with a faster and better transient response. However undesirable results were attained due to lack of well-defined safety features and time delay of the switching rule lead in some of the experiments. \cite{fancher1996comparative} compares Fuzzy logic and H-infinity approaches for ACC. Overall the performances are slow in these designs, and Fuzzy logic design spends a considerable time performing under the desired velocity. Neural networks, and proportional derivative (PD) type control laws were used in \cite{guvencc2006adaptive,breuer1999versatile}.  Model based control for an implementation of an intelligent cruise control was examined in \cite{girard2005intelligent}. This design performs well and the vehicle speeds match the desired values well in performed experiments. Raza and Ioannou in \cite{raza1996vehicle} presented a high-level supervisory control design for vehicle longitudinal control. In \cite{shladover2001nonlinear}, the authors use optimal dynamic back-stepping control to derive the desired acceleration on the supervisory level. However, most supervisory controllers are based on mathematical models rather than real human behavior. Closer to our work, human behavior is modeled through fuzzy controllers or neuro-controllers for spacing adjustments \cite{dermann1995nonlinear,holzmann1997longitudinal}. Very similar to the idea of passivity, Druzhinina et al. \cite{druzhinina2002adaptive} have designed an adaptive cruise control using a Lyapunov function approach. 
	
	Passivity in the field of control design and analysis as a suitable alternative to other nonlinear control techniques has been studied \cite{ortega2001putting,vzefran2001notion}. Passivity conditions for hybrid and switched systems are discussed in \cite{zhao2008dissipativity,manual2002mechanical} where the energy in each mode of the passive system is bounded and that the composite energy of all modes (active or inactive) are also bounded. It is obvious that this will bring about constraints on the design of controller and consequently certain adequate criteria should be met in order for the controller to be passive. However, passivity of the closed loop system can ensure stability. Our experiments are performed in Matlab/Simulink \cite{guide1998mathworks} and CarSim. CarSim is a commercial parameter-based vehicle dynamics modeling that has been used to model the physical world \cite{manual2002mechanical}. 
	
	\section{Background on Passivity and Dissipativity}
	\subsection{Mathematical Preliminaries}
	In continuous time ($\mathbf{T}\in \mathcal{R^+}$), the space of signals of dimension $m$ with finite energy is denoted by $L_2^m$. A continuous time signal $x: T~\rightarrow~ R^m$ is in the general space $H (x\in\ H)$ ) if the signal has finite $L_2^m$- norm,
	
	\begin{align}
	\label{eq:l2norm}
	\|x||_2^2=\int_{0}^{\infty}x^T(t)x(t)d t  <\infty.
	\end{align}
	
	The extended signal space $L^m_{2e}$ can be defined by introducing the truncation operator. The truncation of a continuous signal $x(t)$ to time $T$ is indicated as $x_T(t)$,

	\begin{align}
	x_T(t) = \left\{ \begin{array}{lllll}
	x(t), \quad& \mbox{if ~~$ t\leq~T$} \vspace{0.03in}\\
	0 , \quad& \mbox{if~~ $ t>T$} 
	\end{array} \right.
	\end{align}

	The inner product of signals $y$ and $x$ over the interval $[0,T]$ in continuous time is denoted as 
	
	\begin{align}
		\label{eq:l2norm}
		\ <x,y>_T=\int_{0}^{T}y^T(t)x(t)d t. 
	\end{align}
	
	A continuous time signal $x: T~\rightarrow~ R^m$ is in $L^m_{2e}$ if
	
	\begin{align}
   \label{eq:l2enorm}
   \|x||_T^2=\int_{0}^{T}x^T(t)x(t)d t < \infty \quad \forall T\in T  .
    \end{align}
    
    A system $H$ is a relation on $H_e$. For $u \in H_e$, the symbol $Hu$ denotes to the image of $u$ under $H$. The signal space under consideration is either the standard $L_2^m$ space or the extended $L^m_{2e}$ space. We use $G(s)$ to denote the transfer function for a single-input-single-output linear system.
    
    \subsection{Background on Passivity Theorem}
    
    A system is called passive if it does not generate energy and only stores or dissipates the energy provided to it. This can be shown by an inequality where the energy supplied to the system by the environment, $<Hu, u>_T$, is an upper bound on the loss of stored energy. If we denote $\beta$ as the stored energy available for extraction, we can hold the alternative perspective that the maximum energy that can be extracted from a system, $-<Hu, u>_T$ is bounded above by the stored energy $\beta$.
    
    \begin{definition} %(\cite{Schaft00,Willems72})
	%\label{def:1}
	A relation  $\mathbf{H: u~\rightarrow~ y}$ is said to be \emph{dissipative} with respect to supply rate $w(u,y)$, if
	\begin{align}
	\label{eq:dissipative}
	\int_{t_0}^{t_1}w(u,y)d t \geq~0,
	\end{align}
	for all $t_1\geq t_0$, and all $u\in\mathcal{U}$.\QEDB
	\end{definition}
	
	\begin{definition}%[\cite{Jie07,Lozano00}]
	Suppose the system $\mathbf{H: u~\rightarrow~ y}$ is dissipative. It is said to be
	\begin{itemize}
		\item \emph{passive} if $w(u,y) = u^Ty$;
		\item \emph{input feed-forward passive} (IFP) if there exists a constant $\nu$ so that $w(u,y) = u^Ty-\nu u^Tu$; we call such a $\nu$ an IFP level, denoted as IFP$(\nu)$;
		\item \emph{output feedback passive} (OFP) if there exists a constant $\rho$ so that $w(u,y) = u^Ty-\rho y^Ty$; we call such a $\rho$ an OFP level, denoted as OFP$(\rho)$;
		\item \emph{input-feedforward-output-feedback passive} (IF-OFP) if there exist constants $\delta$ and $\epsilon$ so that $w(u,y) = u^Ty-\delta y^Ty -\epsilon u^T u$; we call such $\delta$ and $\epsilon$ passivity levels, denoted as IF-OFP$(\epsilon,\delta)$;
		\item \emph{finite-gain $\mathcal{L}_2$ stable} if there exists a constant $\gamma\neq 0$ so that $w(u,y) = \gamma^2 u^Tu- y^Ty$, denoted as FGS($\gamma$).
	\end{itemize}
	
	Further, if $\nu>0$, then the system is said to be \emph{input strictly passive} (ISP); if $\rho>0$, then the system is said to be \emph{output strictly passive} (OSP). Similarly, if $\delta>0$ and $\epsilon>0$, then the system is said to be \emph{very strictly passive} (VSP). The largest IFP level $\nu$ is called the \emph{IFP index} and the largest OFP level $\rho$ is called the \emph{OFP index}, respectively. \QEDB
	
	\end{definition}
	
	\begin{remark}
	if either one of two passivity indices is positive, we say that the system has an 'excess of passivity'; and if either of two passivity indices is negative, we say that the system has a 'shortage of passivity.'
	%Thus, the magnitude of the passivity indices quantify `how far a system is from being passive'.
	\end{remark}

%\subsection{Linear Systems}
	\begin{definition}
	The IFP index for a \emph{stable}\footnote{A function $G(s)$ is called \emph{stable} if it is analytic in the closed right half plane of the complex plane, see e.g. \cite{khalil1996nonlinear}.} linear system $G(s)$ is defined as
	
	\begin{align}
	\label{eq:def_freq}
	\nu(G(s))\triangleq\frac{1}{2}\min_{w\in\mathbb{R}}\underline{\lambda}(G(jw)+G^*(jw)),
	\end{align}
	
	where $\underline{\lambda}$ denotes the minimum eigenvalue and $G^*$ denotes the conjugate transpose of $G$.\QEDB
	\end{definition}
	
	\begin{theorem}[\cite{khalil1996nonlinear}]
	Consider the feedback interconnection of two systems $H$ and $G$.
	\begin{enumerate}
		\item If system $H$ and $G$ are passive, then system $\Sigma$ is passive.
		\item If system $H$ and $G$ are output strictly passive (OSP), then system $\Sigma$ is OSP;
		\item If system $H$ is IF-OFP($\epsilon_1,\delta_1$) and system $G$ is IF-OFP($\epsilon_2,\delta_2$), where $\epsilon_1 + \delta_2>0,~ \epsilon_2 + \delta_1>0$, then system $\Sigma$ is finite gain stable (FGS).\QEDB
	\end{enumerate}
	\end{theorem}
	
	It is not necessary for both systems to be passive to guarantee the stability of interconnection. As an example if $H$ is not passive and $\epsilon_1<0$ and $\delta_1<0$ then for the interconnection to be finite gain stable we only require that the system G has the passivity levels: $\epsilon_2>-\delta_1>0$ and $\delta_2>-\epsilon_1>0$.
	
	\section{Passivation Results}
	
	Many different methods in literature such as series, feedback, or parallel interconnections are used to passivate non-passive systems \cite{bao2007process,kelkar1998robust}.These methods are only capable of working with systems that have certain properties such as stability, and minimum-phase property. Hence, these methods are not effective for delay systems that have the form of $G_0(s)e^{-\tau s}$ and are usually approximated by using Pate approximation method. This report is based on our previous works \cite{xia2014passivity,xia2015guaranteeing} where we introduced a new method of passivation that overcomes the aforementioned limitations in existing methods in literature. Consequently, this method can be used to passivate systems with input-output delays.
	
\begin{figure}[!t]
		\centering
		\includegraphics[scale = 0.5]{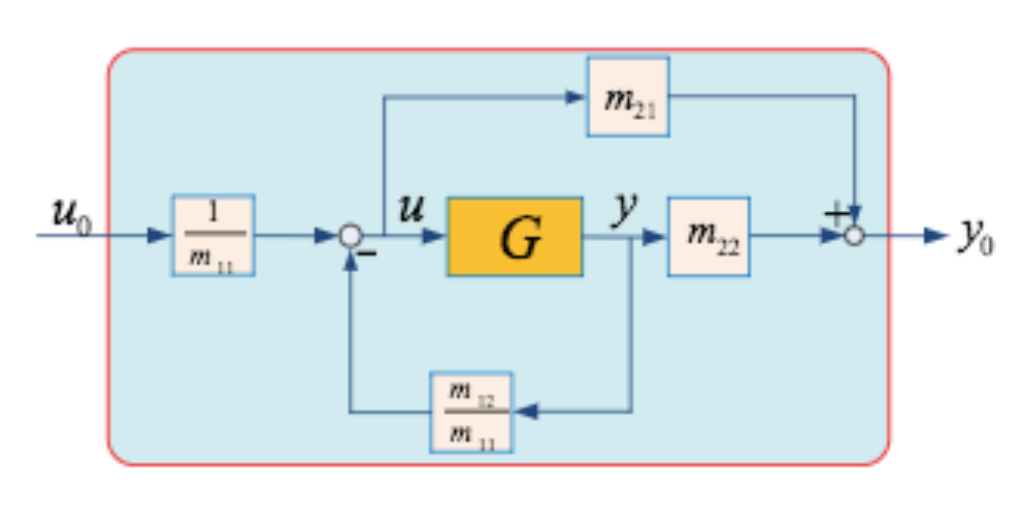}
		\caption{Transformation $M$ given by series, feedforward, and feedback.}
		\label{fig:Ms}
	\end{figure}

	Figure \ref{fig:Ms} shows the passivation method used in our experiment, and more details on this method can be found in \cite{xia2014passivity} and an internal report \cite{xia2014passivityy}.  By considering the non-passive controller $G$ and $m_{11}, m_{12}, m_{21}, m_{22} \in R$ as passivation parameters, we attempt to passivate the system. The interconnection in figure \ref{fig:Ms} can be described as follow:
	
	\begin{align}
	\label{eq:interconnection}
	\begin{bmatrix}
	u_0\\y_0
	\end{bmatrix}=
	M
	\begin{bmatrix}
	u\\y
	\end{bmatrix},
	\end{align}
	
	where the matrix M which is constraint to be invertible is defined as follow:
	
	\begin{align}
	\label{eq:M}
	M=
	\begin{bmatrix}
	m_{11}I&m_{12}I\\m_{21}I&m_{22}I
	\end{bmatrix}.
	\end{align}
	
	In our previous works, we showed that by appropriate selection of passivation parameters we can achieve the desired passivity levels $(\upsilon_0, \rho_0)$ for the system $\Sigma_0: u_0~\rightarrow~y_0$. In this report we focus on passivation using constant gains. In general, it is possible to use transfer functions to replace the constant parameters \cite{xia2014passivityy}.
	
	\begin{theorem}[\cite{xia2014passivity}]
		
	consider a system $G$ which is finite gain stable with gain $\gamma$ and a passivation matrix $M$ as shown in figure \ref{fig:Ms}. Then the system $\Sigma_0: u_0~\rightarrow~y_0$ is
	
		\begin{enumerate}
			\item passive, if $M$ is chosen such that
				\begin{align}
				\label{eq:pm}
			m_{11}=m_{21},\quad m_{22}=-m_{12},\quad m_{11} \geq m_{22}\gamma > 0.
				\end{align}
			\item OSP with OFP level $\rho_0=\frac{1}{2}(\frac{m_{11}}{m_{21}}+\frac{m_{12}}{m_{22}})>0$, if $M$ is chosen such that
				\begin{align}
				\label{eq:opm}
				m_{11}m_{22}>m_{12}m_{21}>0,\quad  m_{21} \geq m_{22}\gamma > 0.
				\end{align}
			\item ISP with IFP level $\nu_0=\frac{1}{2}(\frac{m_{21}}{m_{11}}+\frac{m_{22}}{m_{12}})>0$, if $M$ is chosen such that
		\begin{align}
		\label{eq:opm}
		m_{12}m_{21}>m_{11}m_{22}>0,\quad  m_{11} \geq m_{12}\gamma > 0.
		\end{align}
		\item VSP with passivity levels $\delta_0=\frac{1}{2}\frac{m_{11}}{m_{21}}>0$, and $\epsilon_0=\frac{a}{2}\frac{m_{21}}{m_{11}}>0$ if $M$ is chosen such that
			\begin{align}
			\label{eq:opm}
			m_{11}>0,\quad m_{12}=0, \quad  m_{21} \geq \frac{m_{22}\gamma}{\sqrt{1-a}}>0,    
			\end{align}
				where $0<a<1$ is an arbitrary real number.     
		\end{enumerate}
	
	\end{theorem}
	
		In theorem 2 is $G$ can be nonlinear or linear provided that it is finite gain stable. In our experiment $G$ is a non-passive adaptive cruise controller. 
		
	\section{Co-Simulation Framework}
	
	An extremum-seeking co-simulation based framework was used to select the passivation parameters $m_{11}$, $m_{12}$, $m_{21}$, $m_{22}$ to both passivate the system, and optimize the overall performance. Since the mathematical relation between performance and the parameters is not known, and the relationship is difficult to obtain, we cannot calculate the values directly. As a result we use an extremum-seeking co-simulation based framework to obtain the desired values. Extremum-seeking optimization method is discussed in the next section, and a brief description of co-simulation control structure is given in this section.
	
	Increase in computational capability of computers, coupled with new mathematical programing techniques has produced a particularly promising methodology called co-simulation. The main purpose of this method is to address the control problem rising from the ever-increasing scale and complexity of large-scale systems \cite{baldi2014plug}. Under co-simulation, the system is controlled under realistic conditions including physical constraints and uncommon and possible unpredictable behaviors. Optimization plays a strong role in co-simulation, as most co-simulation frameworks utilize an optimizer to minimize or maximize a cost function which is related to the system performance or control behavior. Under co-simulation, the exact or approximated state-space model of the system is not required, and the control parameters are decided by utilizing the input/output data coming from the system.
    
\begin{figure}[!t]
    		\centering
    		\includegraphics[scale = 0.5]{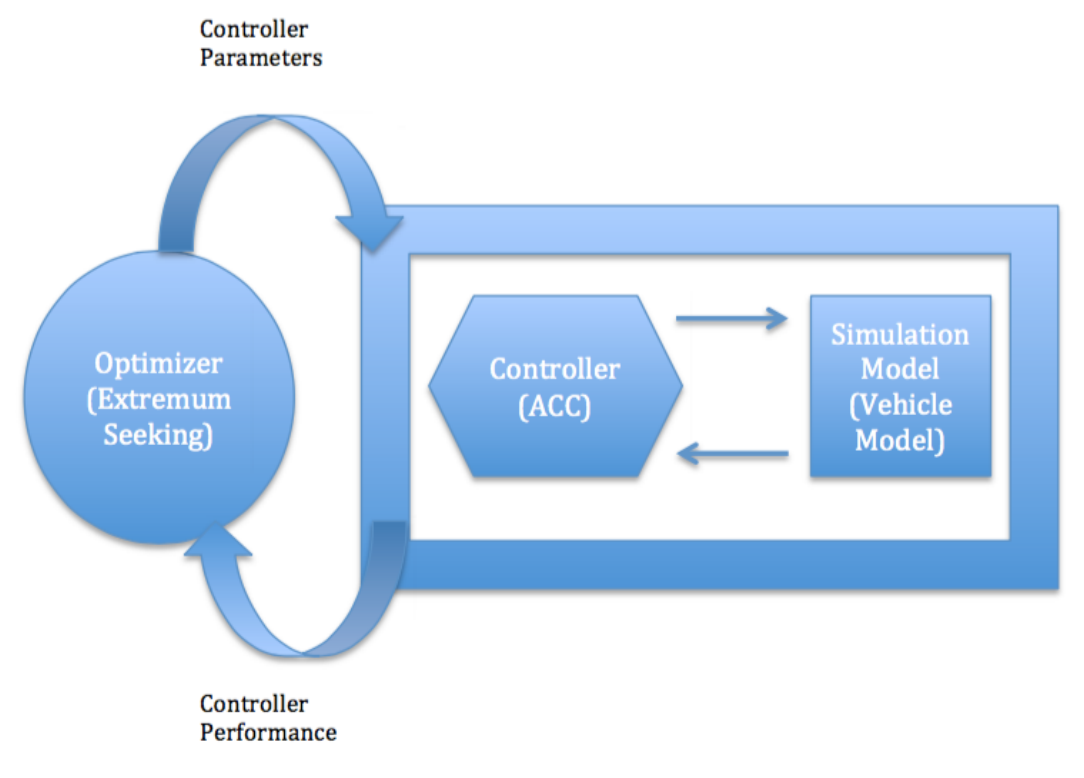}
    		\caption{Extremum-Seeking Based Co-Simulation Framework with Adaptive Cruise Control.}
    		\label{fig:cosim}
    	\end{figure}

    Under co-simulation control structure, a simulation model of the system is connected to a parameterized controller with parameters that are updated by an optimization algorithm. In our experiment, these parameters also passivate the ACC. For each choice of parameter sets, the close loop system is simulated and a measure of the performance is provided to the optimizer \cite{baldi2014plug}. Due to an absence of a well-defined state-space model for the system, the traditional gradient methods employing Jacobian and Hessian matrices are rendered useless for co-simulation algorithms. Consequently, derivative-free optimization methods such as extremum-seeking method play an important part in designing co-simulation controllers. Figure \ref{fig:cosim} shows the basic structure for a co-simulation model used in our experiment.
    
    The main difficulty in co-simulation optimization is the trade-off that exists between allotting computational resources for searching the solution space versus running additional simulation replications for attaining a better estimation of the performance of current solutions \cite{fu2008some}. The general setting for the optimization process used in co-simulation is to find a configuration or design that minimizes the objective function:
    
    	\begin{align}
    	\label{eq:expec}
    	\min\hat{J}=\min E[(J(\theta))]=\min E[L(\theta,\omega)].
    	\end{align}
    	
   where $\theta \in \Theta$ represents the vector of variables to be optimized. In our experiment these variables are $\theta=[m_{11},m_{12},m_{21},m_{22}]$. $J(\theta)$ is the objective function, $\omega$ represents a sample path (simulation replication), and $L$ is the sample performance measure. $\hat{J}$ represents an estimate for expectation of $J(\theta)$. In our experiment the cost function (13) is the tracking error of the form:
   
    	\begin{align}
    	    %	\label{eq:error}
    	\begin{split}
       	\min J=\int_{0}^{T}|v_h-v_{des}|d t \\
        w.r.t \quad m_{11},m_{12},m_{21},m_{22}\\
        subject\quad to: being\quad passive
        \end{split}
    	\end{align}
  
  	where in our adaptive cruise control framework, $v_d$ is the instantaneous desired velocity, and $v_h$ is the instantaneous host velocity. In case of space control, these variables are replaced by instantaneous distance between vehicles, and its desired location, note that there is a new set of $m$'s that passivate the space control:
  	
  	\begin{align}
  	%	\label{eq:error}
  	\begin{split}
  	\min J=\int_{0}^{T}|x_h-x_{des}|d t \\
  	w.r.t \quad m_{11},m_{12},m_{21},m_{22}\\
  	subject\quad to: being\quad passive.
  	\end{split}
  	\end{align}
  	
  	Depending on the distance between the host and lead vehicles the controller switches between the velocity controller and spacing controller. The main goal of velocity control is to achieve and maintain the desired velocity, and the main goal of space control is to decelerate the vehicle, and maintain the safe distance when the distance between two vehicles is less than the desired safe distance e.g., when the lead vehicle brakes or slows down. The switching logic between the controllers depend on the acceleration calculated by two controllers, a negative acceleration leads to a brake command and a positive acceleration will lead to a throttle command \cite{eyisi2013model}.  
  	
	\section{Extremum-Seeking Method of Optimization}
	\subsection{Problem Description}
	
	Due to technological advances, real-time optimization has experienced a resurgence of interest in recent years. The constant increase in complexity of engineering systems, including feedback systems and combined sub-units, has led to many optimization and control challenges. Finding explicit solutions to optimization and control problems for multi-agent, nonlinear and infinite dimensional systems and cyber-physical systems are difficult. Sometimes the presence of competing goals, high dimensionality of the system combined with system's uncertainties adds to this complexity. In many cases if a model-based solution exists, this solution is usually too conservative due to modeling shortcomings. In presence of model uncertainties, the design of a robust input-output strategy to control systems is necessary. High-gain observers \cite{freidovich2008performance} and linear robust controllers \cite{franco2006robust} in combination with the feedback linearization techniques \cite{benosman2014extremum} are some of the existing methods in the literature for handling these situations. Another approach to deal with model uncertainties is using adaptive control methods \cite{krstic1995nonlinear}. The main existing idea behind most real-time adaptive control approaches used for control of linear \cite{aastrom2013adaptive,goodwin2014adaptive} and nonlinear \cite{ioannou1995stable,krstic1995nonlinear} systems only supports the regulation of a known set of points or tracking of a known set of reference trajectories.  A survey of some of these methods is given in \cite{astrom1995theory,aastrom1993automatic}. In most traditional model-based approaches, the model is first adapted using available measurements and then a new numerical optimization is performed on the updated model \cite{zhang2002real}.
	
	However, as mentioned before, in most realistic engineering and scientific problems, the control objective is to optimize a cost function that is usually a function of unknown parameters and variables. Sometimes the explicit form of the objective function is not available to the designer but instead one can measure the value of the output for the objective function by probing the system. In some applications the system needs to operate close to an extremum of a given objective function during its steady-state operation. In other words, a common control goal might be to select a set of desired states that would keep a performance function at its extremum value. In many applications, this set of states or points should be selected such that the system would reach a maximum or minimum point for an uncertain reference-to-output equilibrium map and maintain this behavior. The uncertainty in the input-to-output map makes it necessary for the design to follow an adaptive structure in order to find the set point, which maximizes or minimizes the objective cost function.  
	
	\subsection{History, Background and Related Work}
	
	Two main methods of control used to handle the aforementioned problems are called “extremum-control” and “self-optimization control.” If kept at constant set points, the main task of self-optimizing control is to find a set of controller variables' values that lead the system to near optimal operation \cite{morari1980studies,skogestad2000plantwide}. The main task of extremum-seeking control (also named as external control or peak-seeking control) is to find the operating set points that maximize or minimize an objective function. The idea behind extremum-seeking control was first initiated by Leblanc’s paper in 1922 \cite{leblanc1922electrification}.  Extremum-seeking is a non-model-based real-time optimization approach for many dynamical problems. The first task of extremum-seeking control is to seek an extremum of the performance function, and the second task of this controller is to be able to control (stabilize) the system and drive the performance output to that extremum.
	
	Extremum-seeking control (ESC) is an optimal control approach that deals with situations in which the plant model or the cost function to be optimized is not available but it is assumed that measurements of plant's input and output signals are available.  Extremum-seeking control experienced great popularity long before the theoretical advances and popularity of adaptive linear control in 1980’s. This simple and effective optimization scheme gained further attention in Russia in 1940’s \cite{kazakevich1943technique,kazakevich1944extremum,meerkov1968asymptotic}. In the US, 1950’s and 1960’s witnessed an increase in publications related to extremum seeking control \cite{meerkov1968asymptotic,morosanov1957method,obabkov1967theory,ostrovskii1957extremum,blackman1962extremum,draper1951principles,kazakevich1961extremum}.  Due to lack of a rigorous performance and design analysis of extremum-seeking, this method almost disappeared for several decades until the emergence of a Lyapunov based proof for its stability analysis \cite{luxat1971stability}, with a succeeding resurgence of interest in extremum-seeking for further theoretical developments. In 2000 the first general stability analysis of the extremum-seeking method for stable dynamic systems with unknown output functions was performed by Krstic' and Wang \cite{krstic2000stability}. It has been shown that for any initial condition, one can tune the extremum-seeking controller such that the system parameters converge to a small neighborhood of its local best performance value. In other words the performance of the plant converges to an arbitrary small neighborhood of the optimal steady-state performance \cite{krstic2000performance}. Krstic' and Wang’s work mostly focused on local convergence. Recently, new theories on semi-global stability and non-local convergence of extremum points for arbitrarily large domains have been developed \cite{tan2006non}. Some of the advantages of extremum-seeking control are its robustness against disturbance, and noise, and its relatively less sensitive overall structure. Extremum-seeking is also considered as an optimization approach. Many ideas used in extremum seeking control have been transferred from the field of numerical optimization. Moreover, one can look at extremum-seeking approach, as a kind of constrained optimization problem whose constraints are the system's differential equations. All these qualities make extremum-seeking control a good candidate for automotive systems, where lack of exact physical models is ubiquitous and robustness against noise is crucial. 
	
	Developments in computational power of computers and theoretical advancements in optimization methods such as extremum-seeking have led to an increasing number of applications in industry that use ESC as their main component. Some of these applications are biochemical reactors \cite{guay2004adaptive,bastin2009extremum}, variable cam timing engine operations \cite{popovic2006extremum}, electromechanical valves \cite{peterson2004extremum}. In more recent years, one can name anti-lock braking system controls \cite{drakunov1995abs}, mobile sensor networks \cite{ogren2004cooperative,biyik2008gradient}. \cite{aastrom2013adaptive} and a more recent reference \cite{ariyur2003real} deliver a comprehensive survey of the current literature, history, applications and future of extremum-seeking control. 
	
	\subsection{Main Scheme}
	Nonlinear programming-based extremum seeking method and adaptive control extremum-seeking method shape up the two most dominant algorithms used in extremum-seeking control problems.  In adaptive control methods, a range of adaptive controllers solves the extremum-seeking problem for a large class of systems \cite{aastrom2013adaptive} . The controller makes use of certain excitation (dither) signals, which provide the desired sub-optimal behavior under the hypothesis that the controller parameters are tuned appropriately. Nonlinear programming-based extremum-seeking control method makes use of a constant probing signal that is applied successively to the system to approximate and generate the gradient of the reference-output map (readout map), accordingly classical numerical optimization methods are used to update the system \cite{tee2001solving}. Under gradient-based extremum-seeking method, the controller takes advantage of the gradient information of the readout map to achieve the online optimization goal. This is under the premise that the gradient information will be sufficient to determine an extremum along a convex surface of a performance function. The principle idea behind this method is that an extremum point has a gradient with the magnitude of zero. By analyzing the relationship between the input and output of the cost function, the controller can distinguish whether the current operating point is on the right side or left side of the optimal point, and calculate this distance. As a result, control efforts can be generated to adjust the operating point to converge to the optimal point. Gradient-based extremum seeking control is studied extensively in \cite{banavar2003convergence}. Nonlinear programming based methods are explained in \cite{tee2001solving}. Extremum-seeking control based on sliding mode is detailed in \cite{yu2002extremum,yu2003smooth}. Multi-parametric extremum-seeking adaptive control for nonlinear system with time-varying uncertainties is explored in \cite{benosman2014extremum,benosman2014multi}. Specifically, in \cite{benosman2014extremum}, a robust nonlinear state feedback is designed to make the nonlinear system input to state stable, and then an extremum-seeking scheme is used to overcome the uncertainties. The authors in \cite{benosman2014extremum} improve the performance of a model-based controller by combining it with a robust extremum-seeking methodology. This design seeks to control a linear time-invariant system with structural uncertainties by making use of extremum-seeking’s robustness against uncertainties. In \cite{nevsic2013non} the authors apply a sampling optimization method called the discrete-time Shubert algorithm to continuous-time plants, and explore the region of convergence for the system and show that semi-global convergence can be achieved in the presence of local extrema.  \cite{khong2013multidimensional} expands this idea to multi-dimensional domains and proves semi-global convergence for the systems based on a periodic sampled-data control law. Additionally the authors explore the robustness for their design.	
	
	In the signal perturbation methods, the process is slightly distorted by some kind of perturbation signals (periodic sinusoidal function of time) to extract the gradient information of the input-output process \cite{krstic2000performance}. In comparison to switching methods \cite{pan2003stability}, which are only concerned with the sign, the perturbation-based extremum-seeking algorithms achieve a better performance by using the gradient sign and magnitude information to set the adaptation speed. The perturbation strategy can save energy by adoptively seeking for a magnitude level of control gain large enough to neutralize all disturbance and uncertainties. Additionally, this method of control is continuous and avoids high feedback control gains such that the system would not experience any chattering phenomenon. 
	
	The cost function is usually assumed to provide the designer with output values based on input signals. The composition of the plant and the cost function are considered as one lumped plant \cite{wang1999optimizing}. It is assumed that an unknown static map can give the steady-state input-output relation of lumped system for fixed inputs. This map is assumed to include an extremum for the parameter values being optimized. Accordingly, the job of the controller is to find this extremum point, based on the assumption that by finding the extremum point of the static map, the steady-state performance of the plant is optimized consequently. This static relationship is justified if the dither frequency is small enough and the time between changes in the optimal reference is sufficiently long. Small dither frequency causes very slow convergence rate, which is acceptable for fast systems, however becomes problematic for biochemical systems with time constants that are in hours or days. It is beneficial to mention that in some cases, the plant is defined as a cascade system of a nonlinear static map and a linear dynamic system - this is known in literature as Hammerstein and Wiener models \cite{nmark1995adaptive}.
	
	In 1922, Leblanc proposed a new method for maintaining the maximum power capacity in the transmission lines. The intuitive explanation behind this idea was simple. Starting with some initial condition , the system's input is perturbed in the form:
	
		\begin{align}
		%	\label{eq:lblank1}
       x(t)=x(0)+\cos(\omega t)
		\end{align}
		which results in a perturbation of some unknown output function:
		\begin{align}
			%	\label{eq:lblank2}
			F(x(t))=F(x(0)+\cos(\omega t))
		\end{align}
		
\begin{figure}[!t]
					\centering
					\includegraphics[scale = 0.65]{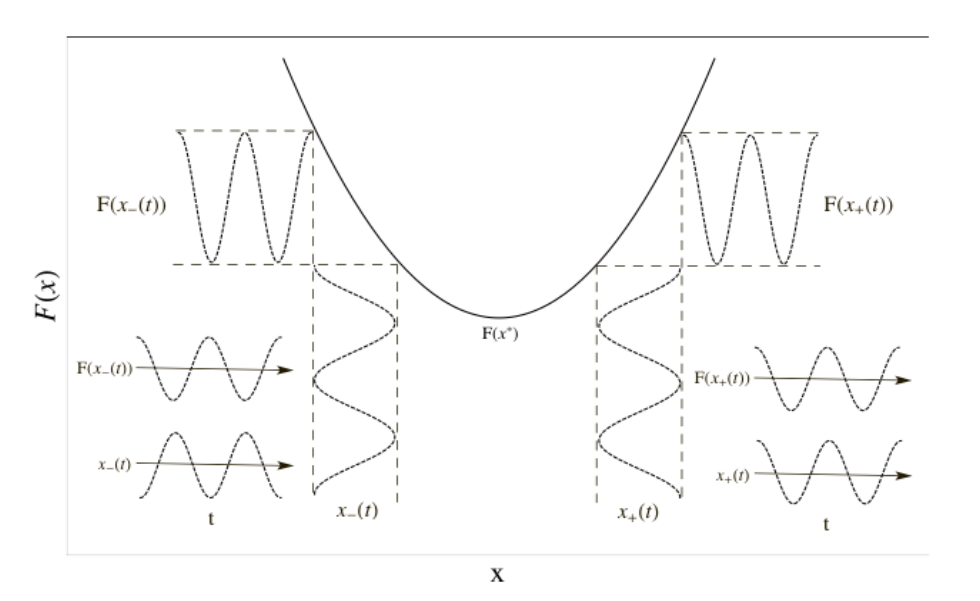}
					\caption{The relationship between perturbation signal, and the extremum point \cite{scheinker2012extremum}.}
					\label{fig:LB}
		\end{figure}

		If the system has an extremum point near this initial point then   will be in or out of phase with the perturbing term  , as depicted in figure \ref{fig:LB}. The phase of the output function, relative to the perturbing signal, gives the controller an idea of how to change the parameter in order to approach the extremum point. Implementing this design in a feedback loop system results in an adaptive automated control system, which leads the system into its minimum or maximum.
			
\begin{figure}[!t]

			\begin{subfigure}{\columnwidth}
					\centering
					\includegraphics[scale = 0.65]{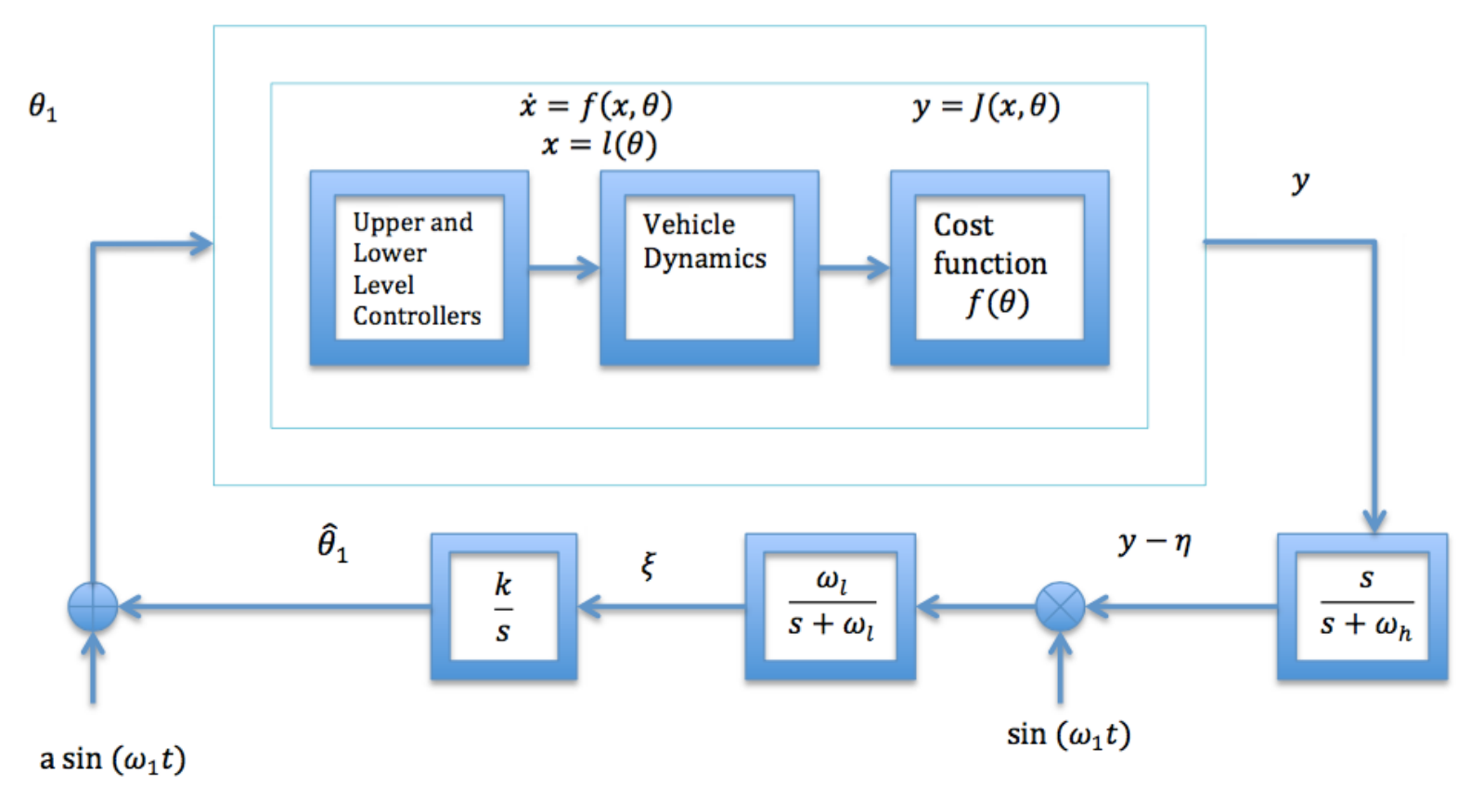}
					\caption{Adaptive Cruise Control and Extremum Seeking Control Design for one parameter.}
					\label{fig:ESschemeA}
			\end{subfigure}
					\hfill	  
			\begin{subfigure}{\columnwidth}
					\centering
					\includegraphics[scale = 0.65]{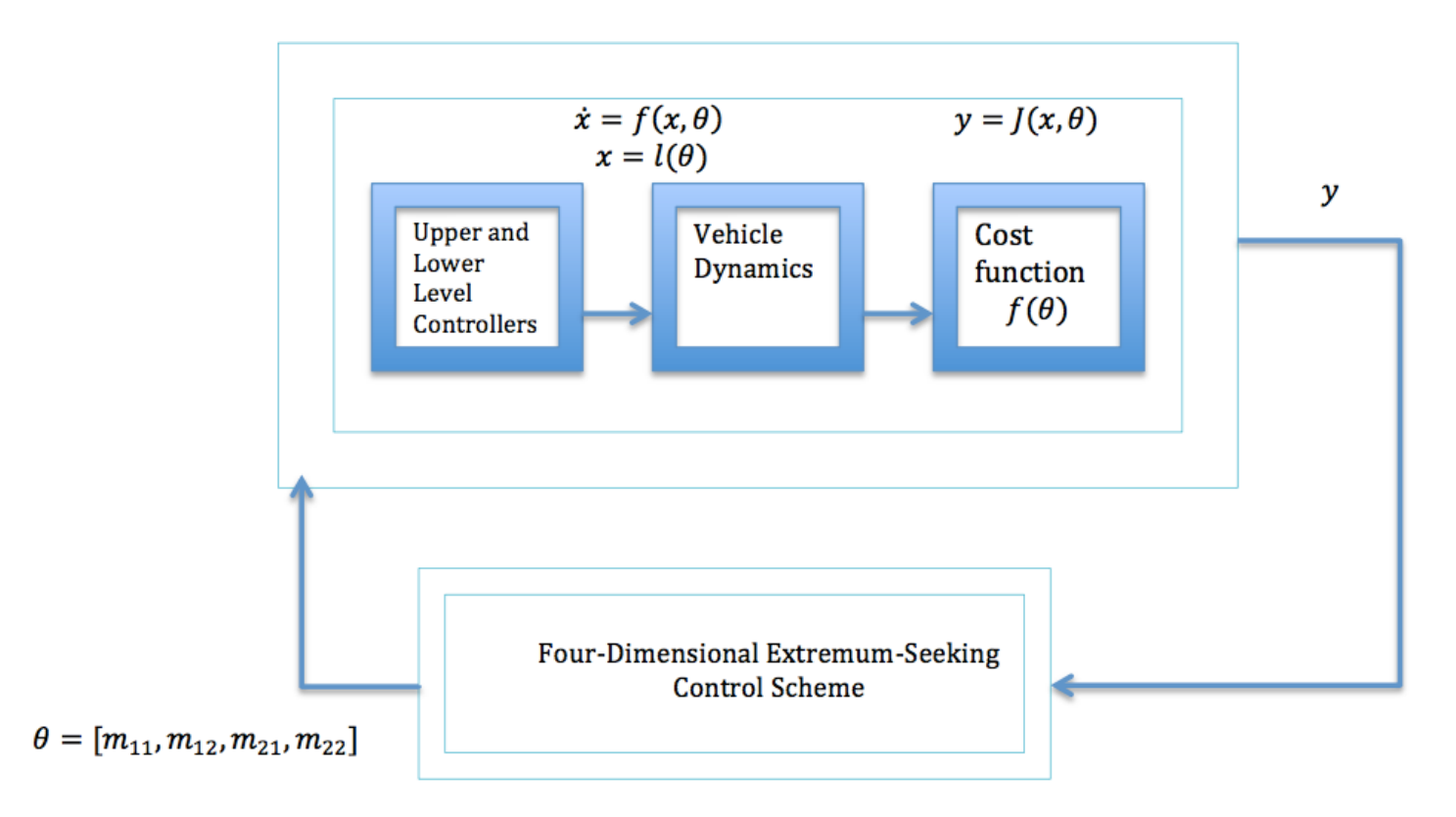}
					\caption{Adaptive Cruise Control and Extremum Seeking Control Design for four parameters.}
					\label{fig:ESschemeB}
			\end{subfigure}
					\caption{Extremum Seeking Control Design.}
			\end{figure}

		Figure \ref{fig:ESschemeA} shows the design used to optimize the performance of our adaptive cruise control for one of the parameters, and figure \ref{fig:ESschemeB} shows this design for optimizing all four parameters. We assume a quadratic structure for the system’s cost function:
	
			\begin{align}
					%	\label{eq:cost}
                    J(\theta)=J^*+\frac{J^{''}}{2}(\theta-\theta^*)^2,
			\end{align}
		
			where $J^{''}>0$ . We know that any $C^2$   function can be expressed in the form (18) with its Taylor series expansion and under the condition that the gradient of $J(\theta)$ is equal to zero for the extremum point. In cases, where $J^{''}<0$ , we can easily replace $k>0$  with $-k$ . The role of the additive probing term $asin(\omega t)$ and multiplicative term of $sin(\omega t)$   along with effects of the washout high pass filter, and the low pass filter is to estimate the derivative of $J$, which is then fed into the integrator in order to perform the classical gradient-based optimization for input-output map with step size $k$ . 
			
			In figure \ref{fig:ESschemeA}  we have, the following expression for estimating the error:
			
			\begin{align}
			%	\label{eq:1}
				\widetilde{\theta_1}=\theta_1^*- \widehat{\theta_1}\nonumber,
			\end{align}
		    so we have:
		    	\begin{align}
		    	%	\label{eq:2}
		    	\theta_1-\theta_1^*=asin(\omega t)- \widetilde{\theta_1}\nonumber,
		    	\end{align}
		    
		    substituting this term into (18), expanding the parenthesis, and using the trigonometric identity $2sin^2(\omega t)=1-cos(2\omega t)$  we get the following:
		    
		      	\begin{align}
		      	%	\label{eq:3}
		      	J(\theta_1)=J^*+\frac{a^2J^{''}}{4}+\frac{J^{''}}{2}\widetilde{\theta_1^2}-aJ^{''}\widetilde{\theta_1}sin(\omega t)+\frac{a^2J^{''}}{4} cos(2\omega t)\nonumber,
		      	\end{align}
		    after the high-pass washout filter, the DC component is removed and signal is demodulated by multiplication with $sin(\omega t)$ resulting in:
		      	\begin{align}
		      	%	\label{eq:4}
		      	\Psi\approx \frac{J^{''}}{2}\widetilde{\theta_1^2}sin(\omega t)-aJ^{''}\widetilde{\theta_1}sin^2(\omega t)+\frac{a^2J^{''}}{4} cos(2\omega t)sin(\omega t)\nonumber,
		      	\end{align}	
		    using both the previous trigonometric identity, and $2cos(\omega t)sin(\omega t)=sin(3\omega t)-sin(\omega t)$ we have:
		    
		    \begin{align}
		    %	\label{eq:5}
		    \Psi\approx -\frac{aJ^{''}}{2}\widetilde{\theta_1^2}+\frac{aJ^{''}}{2}\widetilde{\theta_1}cos(2\omega t) +\frac{J^{''}}{2}\widetilde{\theta_1^2}sin(\omega t)+\frac{a^2J^{''}}{8} (sin(\omega t)-sin(3\omega t))\nonumber,
		    \end{align}	
		    
		    after the low pass filter, all high frequency terms disappear. Given the fact that:
		     \begin{align}
		     %	\label{eq:6}
		     \widetilde{\theta_1}^{'}= -\widehat{\theta_1}^{'} \nonumber,
		     \end{align}	
		     
		     we have: 
		      \begin{align}
		      %	\label{eq:7}
		      \widetilde{\theta_1} \approx \frac{k}{s}[\frac{aJ^{''}}{2}\widetilde{\theta_1}]\nonumber,
		      \end{align}
		      
		      and:
		      
		        \begin{align}
		        %	\label{eq:8}
		        \widetilde{\theta_1}^{'} \approx \frac{-kaJ^{''}}{2}\widetilde{\theta_1}]\nonumber,
		        \end{align}
		        
		     since $kJ^{''}>0$, the design is stable, and we can conclude that $\widetilde{\theta_1}\rightarrow0$ or in other words $\theta_1$ converges to a close neighborhood of $\theta_1^*$.  
		      	
			Another interesting aspect of this extremum-seeking scheme is that the terms $sin(\omega t)$ and $asin(\omega t)$  mimic the job of amplitude modulation (AM) in communications systems. One advantage of extermum-seeking method compared to other optimization methods such as Hooke and Jeeve’s is that ESC performs an online optimization in the sense that it only requires one experiment per iteration to estimate the gradient. This renders ESC systems a much faster performance rate compared to many other optimization schemes.
			
			It is necessary to mention that the design parameter $\omega$ is qualitatively large relative to $k, a, \omega_h$  and that $\omega_l$ should be smaller than $\omega$. Additionally the smaller $a$, the smaller the residual error at the minimum is, but also the larger is the possibility of getting stuck at a local minimum. The gain parameter $k$ jointly with $a$  determine the speed of convergence. The perturbation frequency $\omega$  controls the separation between the time scale of the estimation process performed by the integrator and that of the gradient estimation process performed by the additive and multiplicative perturbation. The higher the frequency $\omega$ the cleaner the estimate of the gradient is.
			
			The average theory for nonlinear non-autonomous systems can be used to examine the stability and convergence of the extremum-seeking systems. In a more qualitative sense, the main idea behind averaging method is that a time-varying periodic system can be approximated by the time-invariant system derived from integrating the original system over a single period. Consequently, if the resulting time invariant system is stable around an equilibrium point, the time-varying periodic system converges to a periodic orbit around this equilibrium point. For more details on the average theory, and its application for stability and convergence of nonlinear systems, one can refer to chapter 9, and more specifically theorem 9.3 in \cite{khalil1996nonlinear}. A detailed proof for convergence and stability of this scheme is given in \cite{krstic2000stability}. The solutions of the closed-loop system converge to a small neighborhood of the extremum of the equilibrium map. The size of this neighborhood is proportional to the adaptation gain, the amplitude and frequency of the periodic signal used in the design. \cite{nesic2009extremum} explains the trade-off existing between the convergence rate and the size of the domain of attraction. Namely, if the domain of attraction is larger, the slower the convergence of the algorithm will be. In \cite{chioua2007dependence} it was shown that the neighborhood of convergence is proportional to the square of the dither frequency. The implication of this is that even when the amplitude of excitation goes to zero, the size of this neighborhood will not go to zero. The following sub-section explains this design and its stability and convergence properties in more details.
			
	 	\section{Simulation Results in CarSim and Matlab/Simulink}		
	 	
	 	In most vehicle control systems, adaptive cruise control in host vehicle following another vehicle adjusts both speed and inter-vehicular distance with respect to the leading vehicle or a moving obstacle for safety spacing purposes. Uncertainties both from surrounding environment and the agents themselves make such a problem challenging. In recent years, the problem of adaptive cruise control has been addressed in different forms such as PID and pole placement control \cite{zhang1999using}, and sliding mode control \cite{lu2002acc}. In our design we combine passivity and extremum-seeking control together to overcome the uncertainties and reach an acceptable performance index. Passivity and extremum-seeking control enable the host vehicle to adaptively seek for a magnitude level of throttle or braking efforts that are necessary to reject disturbance even though the variations of disturbances are unknown. Furthermore, the control input is continuous and thus it will not incur chattering in operation.
	 	
	 	For accurate tracking problems, such as adaptive cruise control, it is usually expected to have enough information about the system's model and reference signal and the system's structure. Given the uncertain environment on the roads, knowledge on system's model may not always be sufficiently available. Hence, extremum-seeking strategy is a good substitution for some class of adaptive control methods. The extremum-seeking scheme explained in previous section can be expanded to multi-parameter cases \cite{ariyur2003real,khong2013multidimensional}. For our four dimensional system, where we optimize the passivation parameters $m_{11}$, $m_{12}$, $m_{21}$, $m_{22}$, the assumption 1 for our design should hold \cite{ariyur2003real}, table 1 shows the design parameters chosen for our first set of experiments. Both frequencies picked for high pass and low pass filters ($\omega_l=1$ and $\omega_h=3$) are smaller than dither frequencies, figure 5 shows the simulation results for $m_{11}$, $m_{12}$, $m_{21}$, $m_{22}$ for this experiment. This experiment lasted for 120 seconds and the control objective was to achieve an asymptotic tracking of the desired velocity.
	
		\begin{assumption}
			\label{as:1}
			$\omega_p+\omega_q \neq \omega_r$ for any $p, q, r= 1, 2 ,3 ,4$
		\end{assumption}	
		
    	\begin{table}[]
		\centering
		\caption{Design Parameters}
		\label{table1}
		\begin{tabular}{|c|c|c|}
            \hline
			 & $\omega$  & $a$ \\  \hline 
			$m_{11}$  & $5+\frac{\pi}{2}$ & $1$ \\ \hline
			$m_{12}$  & $14.14$     & $0.4$\\ \hline
			$m_{21}$  & $7.7$       & $1$ \\ \hline
			$m_{22}$ & $11$        & $0.7$\\\hline
			
		\end{tabular}
    	\end{table}
    	
\begin{figure}[!t]	
			\begin{subfigure}{\columnwidth}
						\centering
						\includegraphics[scale = 0.65]{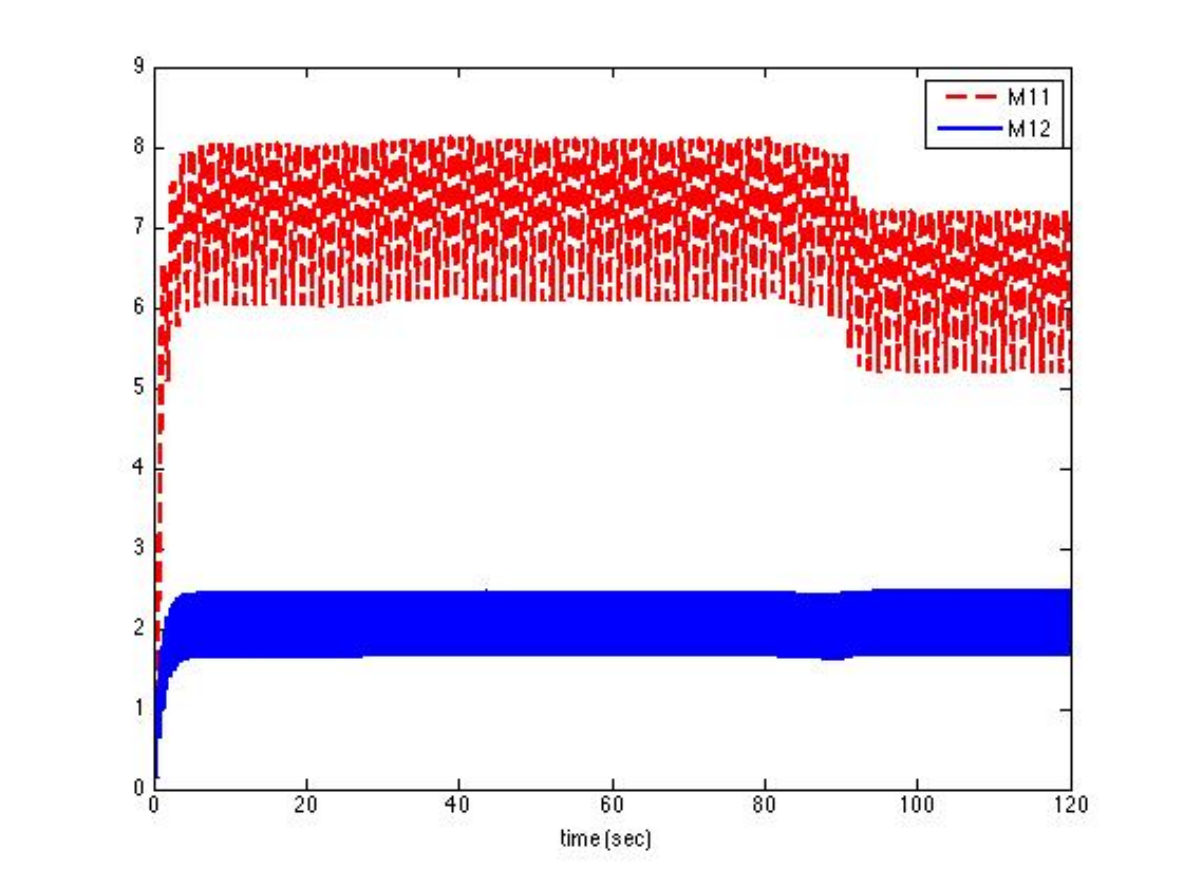}
						\label{fig:firstexpms1}
			\end{subfigure}
					\hfill	  
				\begin{subfigure}{\columnwidth}
						\centering
						\includegraphics[scale = 0.65]{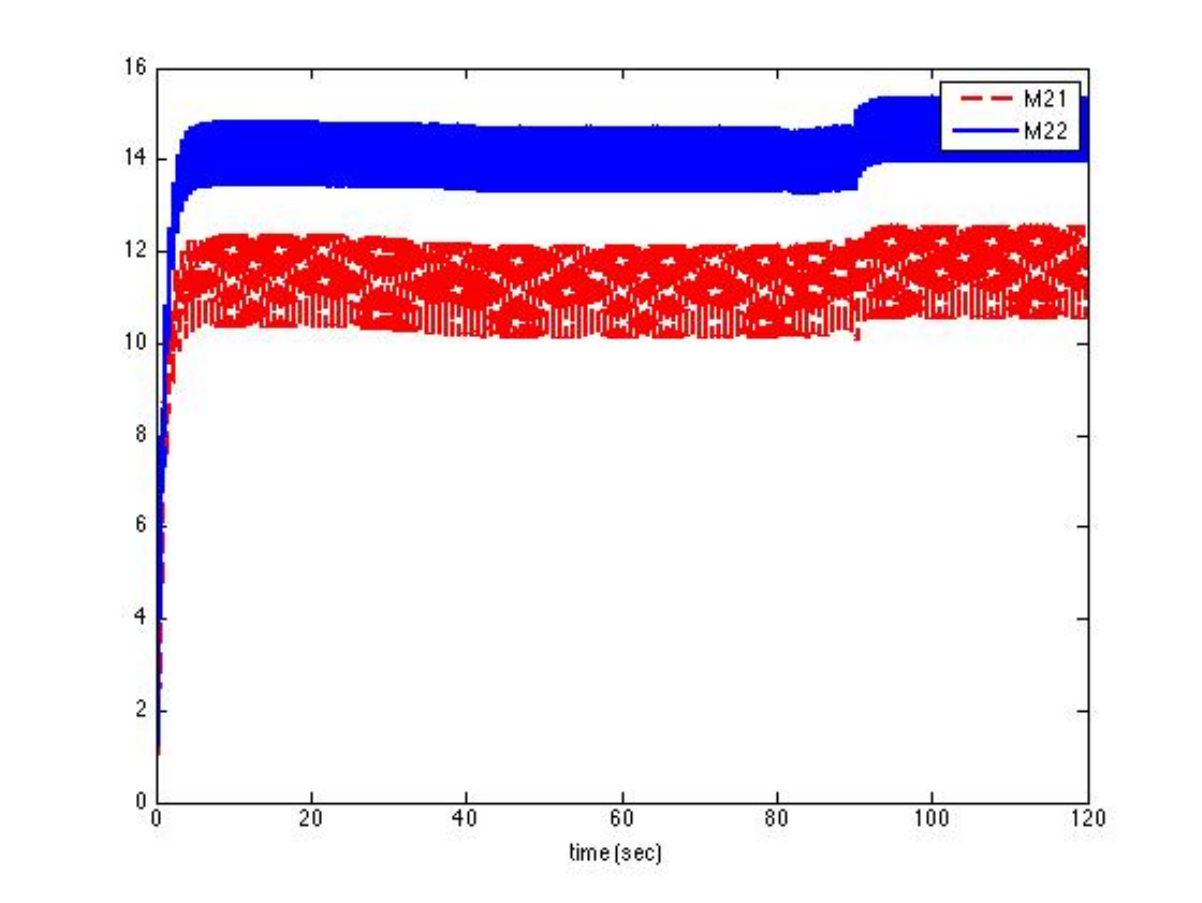}
						\label{fig:firstexpms2}
					\end{subfigure}
					\caption{Optimized Parameters over 120s in First
					 Experiment.}
					\label{fig:firstexpmss}
				\end{figure}

\begin{figure}[!t]
					\centering
					\includegraphics[scale = 0.65]{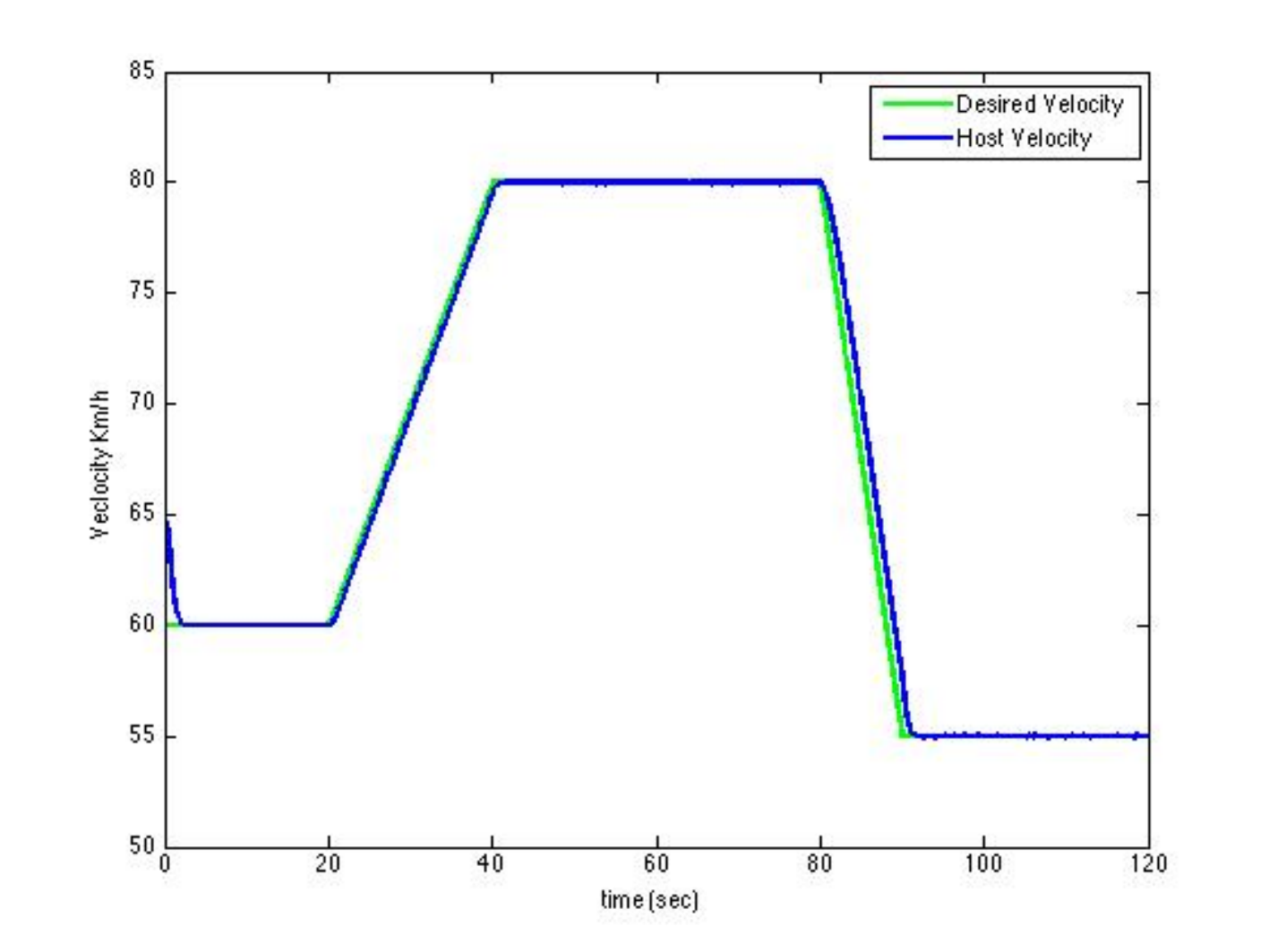}
					\caption{Host Velocity vs. Desired Velocity (120s)     .}
					\label{fig:ACCone}
				\end{figure}

    	The reason for large oscillations in figure \ref{fig:firstexpmss} is the relatively small magnitude of $a$ 's. In general, the magnitude of $a$ is inversely proportional to the amount of oscillation witnessed in estimated parameters. Figure \ref{fig:ACCone}  demonstrates the desired and host velocities for this experiment. As mentioned before the magnitude of $\omega$ inversely affects the rate of convergences, meaning the smaller $\omega$ is, the faster the rate of convergence. The magnitude of $a$ also affects the rate of convergence. The smaller $a$, the larger the domain of attraction, and as a result this means that the system will converge slower.
		
\begin{figure}[!t]	
						\begin{subfigure}{\columnwidth}
							\centering
							\includegraphics[scale = 0.65]{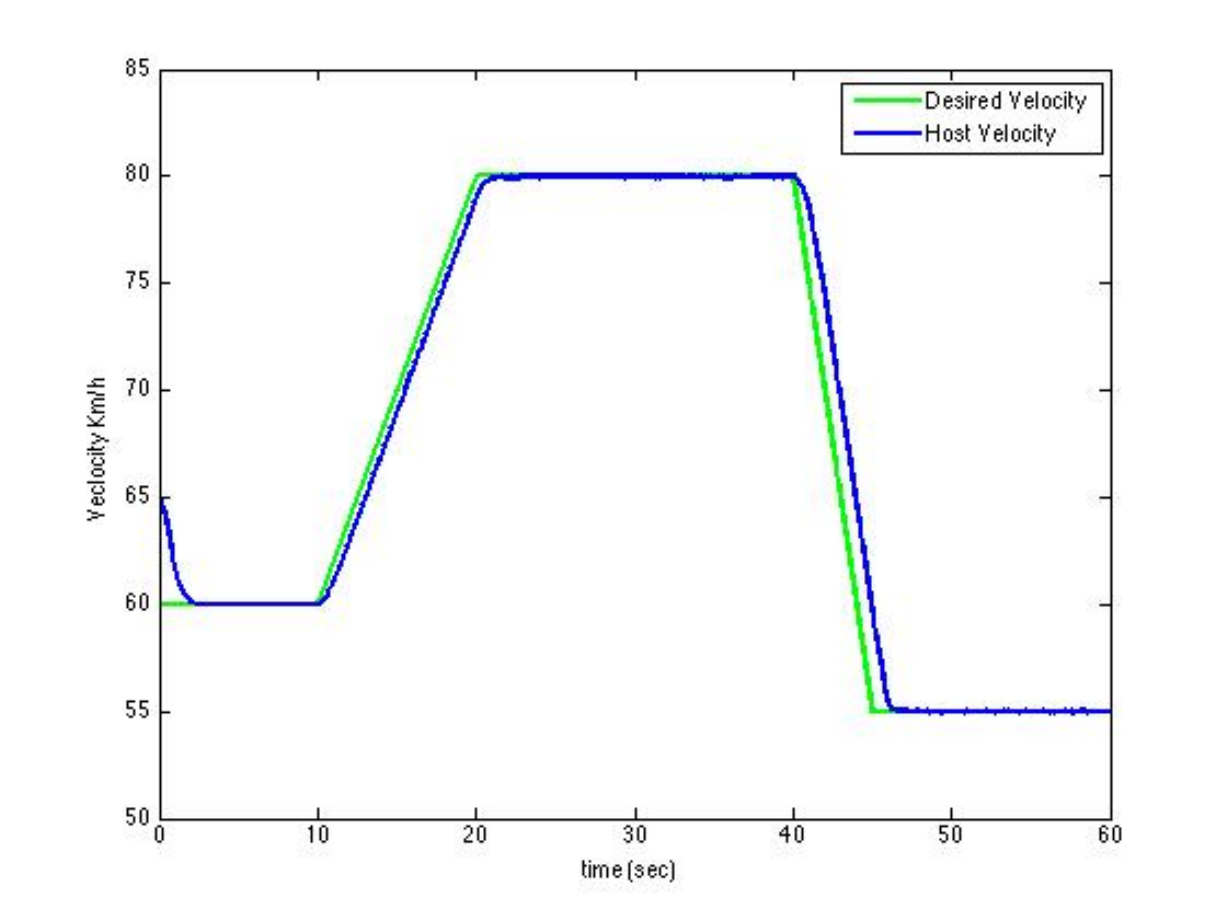}
							\label{fig:ACConeshort}
						\end{subfigure}
						\hfill	  
						\begin{subfigure}{\columnwidth}
							\centering
							\includegraphics[scale = 0.65]{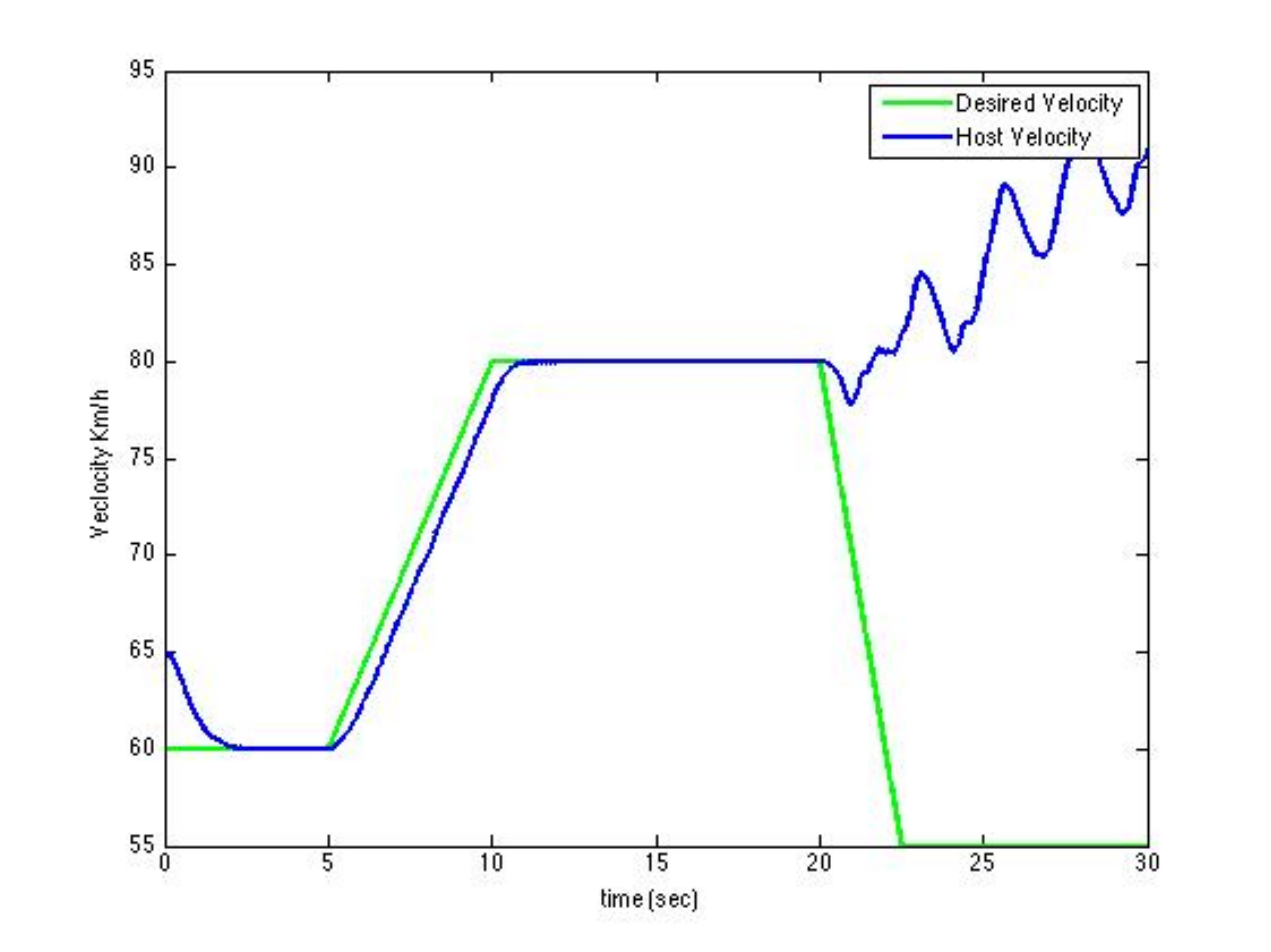}
							\label{fig:ACConeshort2}
						\end{subfigure}
						\caption{Results for 30s long, and 60s long experiments.}
						\label{fig:ACConeshortt}
					\end{figure}

		Figure \ref{fig:ACConeshortt} shows the results for two other experiments in which one can witness the effects of design parameters. Incorporating the same design parameters in the next two experiments: in first experiment that lasted for 60 seconds, one can see that the error has slightly increased compared to figure \ref{fig:ACCone}, given the fact that the shorter time requires a faster rate of convergence. Similarly in the second experiment, since the magnitudes of the set of dither frequencies are not small enough and $a$ is too small to provide fast convergence rate for the 30-second long experiment, the host velocity cannot track the desired velocity correctly and diverges. 
		
\begin{figure}[!t]	
							\begin{subfigure}{\columnwidth}
								\centering
								\includegraphics[scale = 0.65]{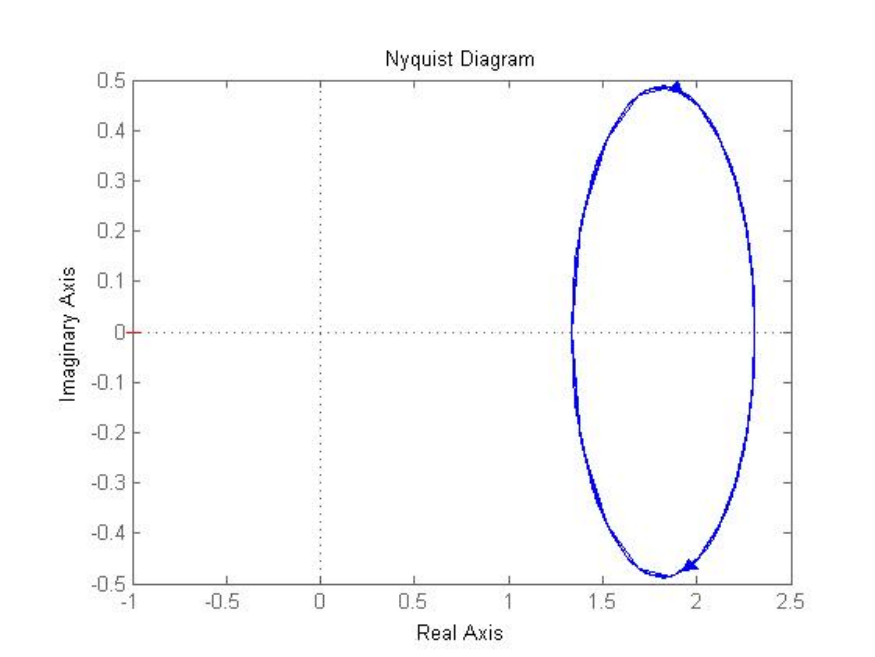}
								\label{fig:Nyqfirst}
							\end{subfigure}
							\hfill	  
							\begin{subfigure}{\columnwidth}
								\centering
								\includegraphics[scale = 0.65]{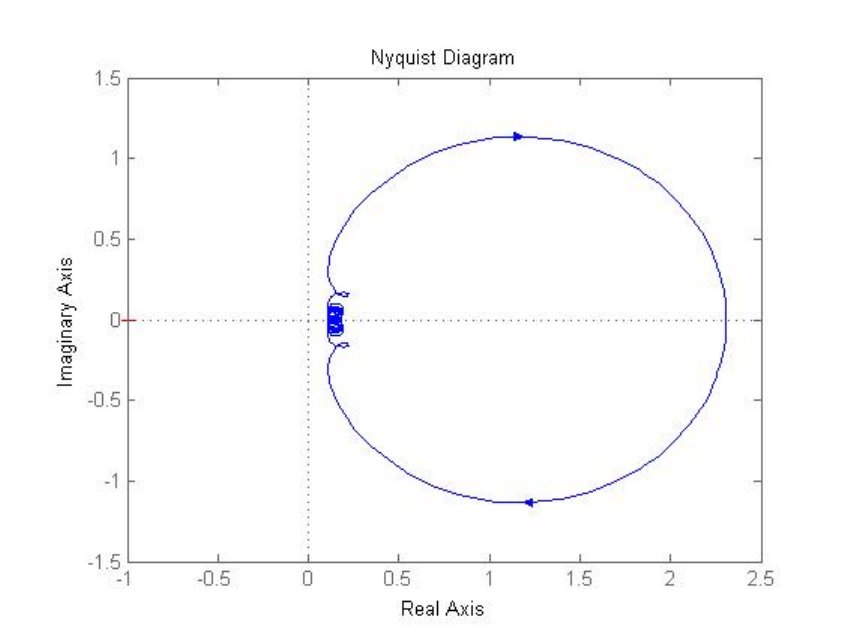}
								\label{fig:Nyqfirst2}
							\end{subfigure}
							\caption{Nyquist and Inverse Nyquist Plots for the passivated system in figure 5, and 6.}
							\label{fig:Nyqfirstt}
						\end{figure}

		Figure  \ref{fig:Nyqfirstt} shows the Nyquist and inverse Nyquist plots that lie on the right-half place for the passivated system. The desired values reached for parameters at the end of experiment are given in table 3. Our controller has the gain of $0.5$, as a result condition (10) holds for the passivated system. The overall design is output passive with output passivity level of $\rho_0=0.385$. 
		
\begin{figure}[!t]
								\centering
								\includegraphics[scale = 0.65]{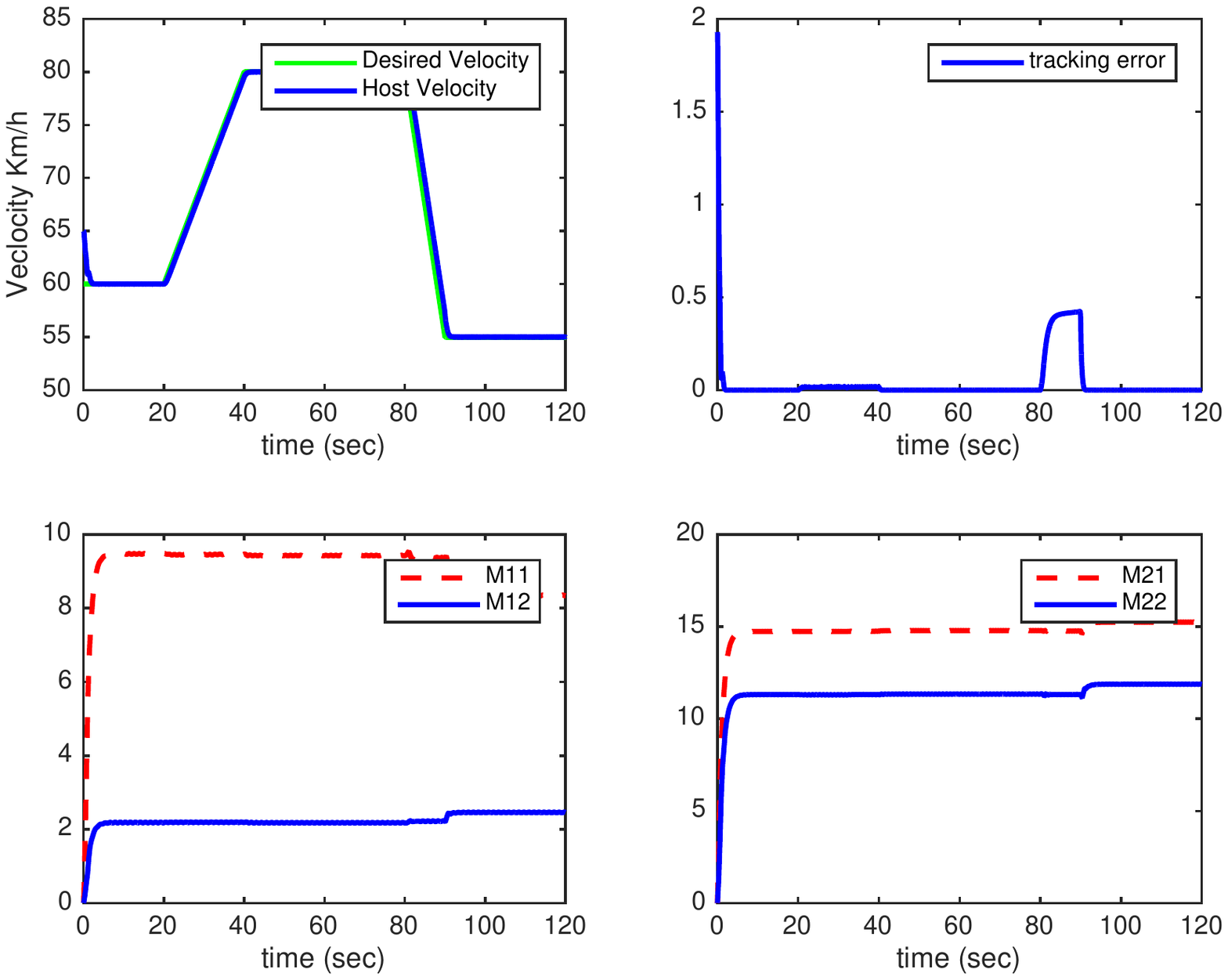}
								\caption{Simulation Results for Second Experiment.}
								\label{fig:secondexp}
							\end{figure}

		Table 2 shows a new set of values used for our second set of experiments, where we achieve a steadier tracking performance. The values for parameters that are being optimized are less oscillatory, and achieve constant values toward the end (figure \ref{fig:secondexp}).
		
		\begin{table}[]
				\centering
				\caption{Design Parameters}
				\label{table2}
		\begin{tabular}{|c|c|c|}
					\hline
					& $\omega$  & $a$ \\  \hline 
					$m_{11}$  & $5+\frac{\pi}{2}$ & $0.01$ \\ \hline
					$m_{12}$  & $5$     & $0.01$\\ \hline
					$m_{21}$  & $10.14$       & $0.01$ \\ \hline
					$m_{22}$ & $7$        & $0.01$\\\hline
					
		\end{tabular}
			\end{table}
		
				\begin{table}[]
					\centering
					\caption{Final values for Parameters}
					\label{table3}
					\begin{tabular}{|c|c|}
						\hline
						$m_{11}$  & $8.6$ \\ \hline
						$m_{12}$  & $3.47$\\ \hline
						$m_{21}$  & $16.51$ \\ \hline
						$m_{22}$ & $13.71$\\\hline
						
					\end{tabular}
				\end{table}
						
\begin{figure}[!t]
				\centering
				\includegraphics[scale = 0.65]{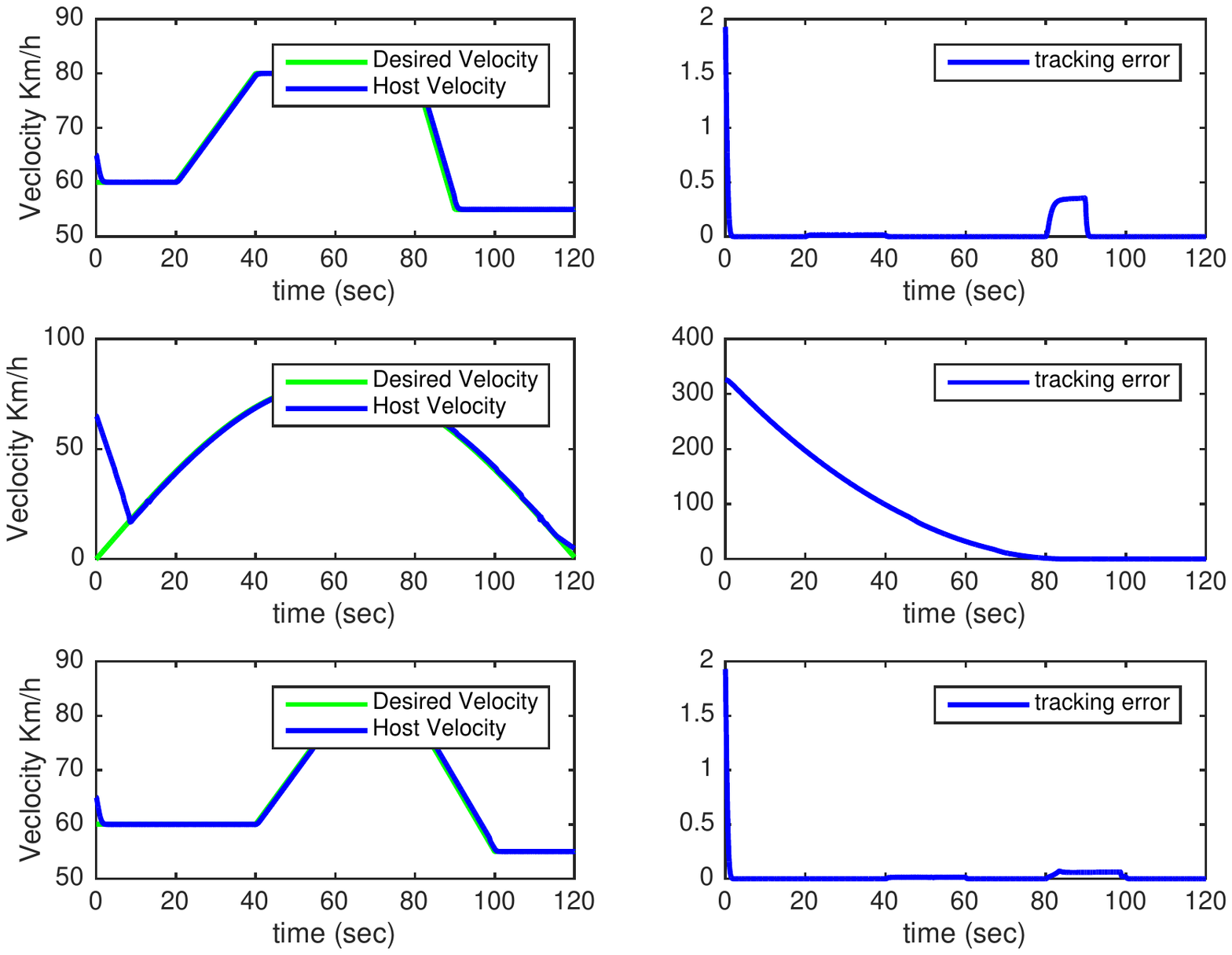}
				\caption{Performance Results for Constant Optimized Control Parameters.}
				\label{fig:diffspeedprofile}
			\end{figure}

	In our next experiment, first we set the control parameters to desired values reached in previous experiment (given in Table 3), then we seek to not optimize the parameters but use parametersthem to examine the controller performance for different speed profiles (figure \ref{fig:diffspeedprofile}). The good performance shown in figure \ref{fig:diffspeedprofile} shows that one does not necessary need to know the shape of the tracking target in advance to select controller parameters accordingly. In other words, the optimized  for a certain experiment  can still improve the behavior of system for different speed trajectories.
	
\begin{figure}[!t]
				\centering
				\includegraphics[scale = 0.65]{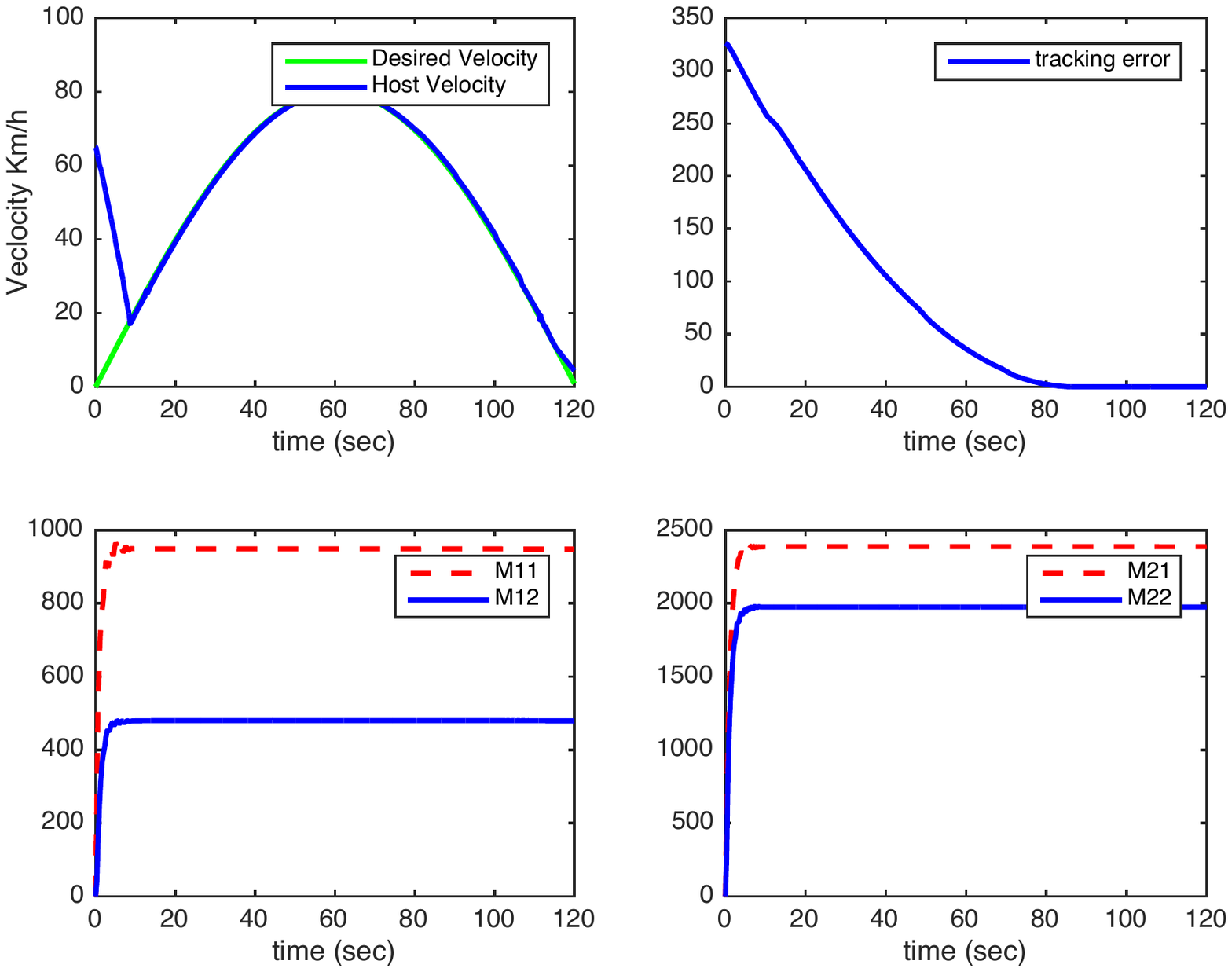}
				\caption{Simulation Results for Sinusoidal Speed Profile.}
				\label{fig:sin}	
			\end{figure}

\begin{figure}[!t]	
					\begin{subfigure}{\columnwidth}
						\centering
						\includegraphics[scale = 0.65]{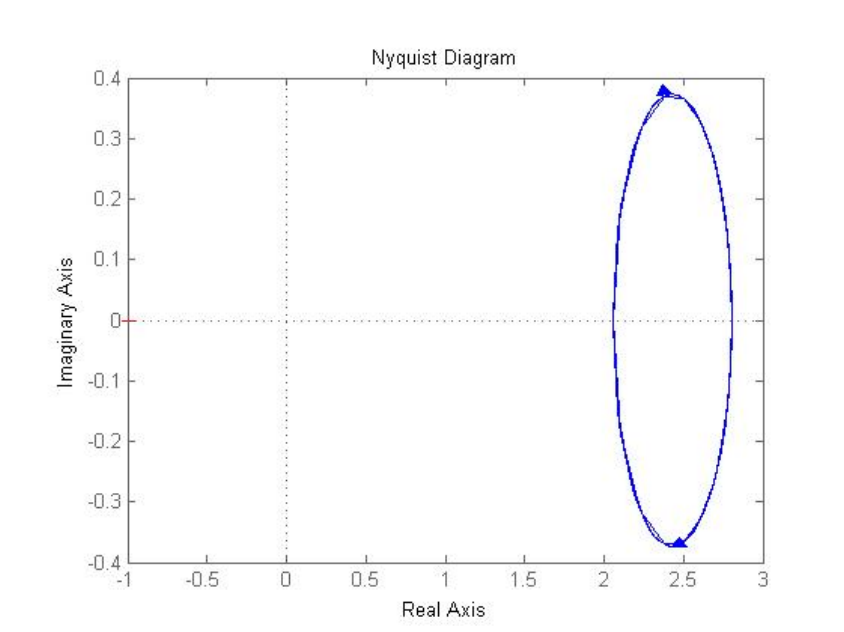}
						\label{fig:NyqSec}
					\end{subfigure}
					\hfill	  
					\begin{subfigure}{\columnwidth}
						\centering
						\includegraphics[scale = 0.65]{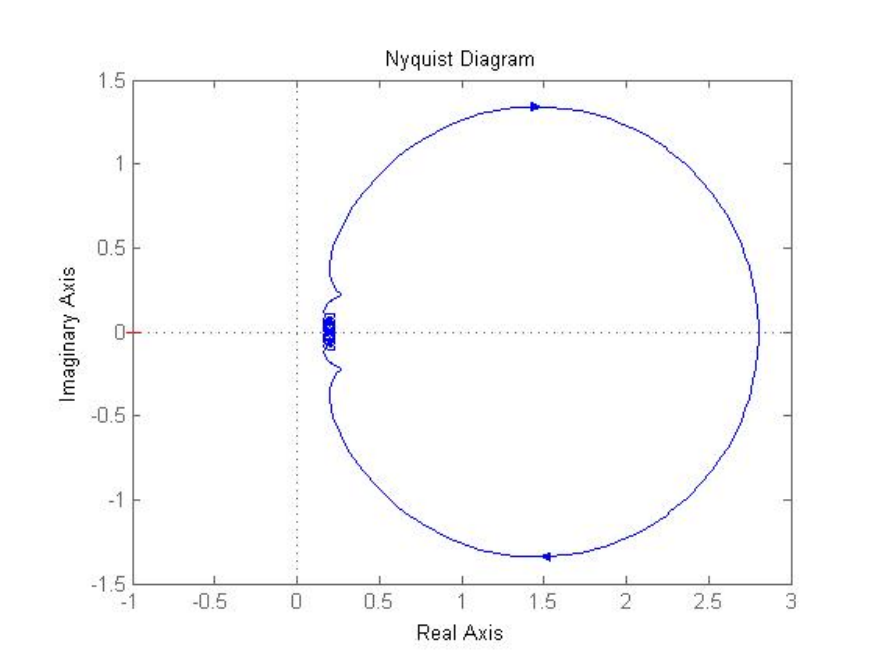}
						\label{fig:NyqSec2}
					\end{subfigure}
					\caption{Nyquist and Inverse Nyquist Plots for the passivated system in figure 11.}
					\label{fig:Nyqsecc}
				\end{figure}

\begin{figure}[!t]
		\centering
		\includegraphics[scale = 0.65]{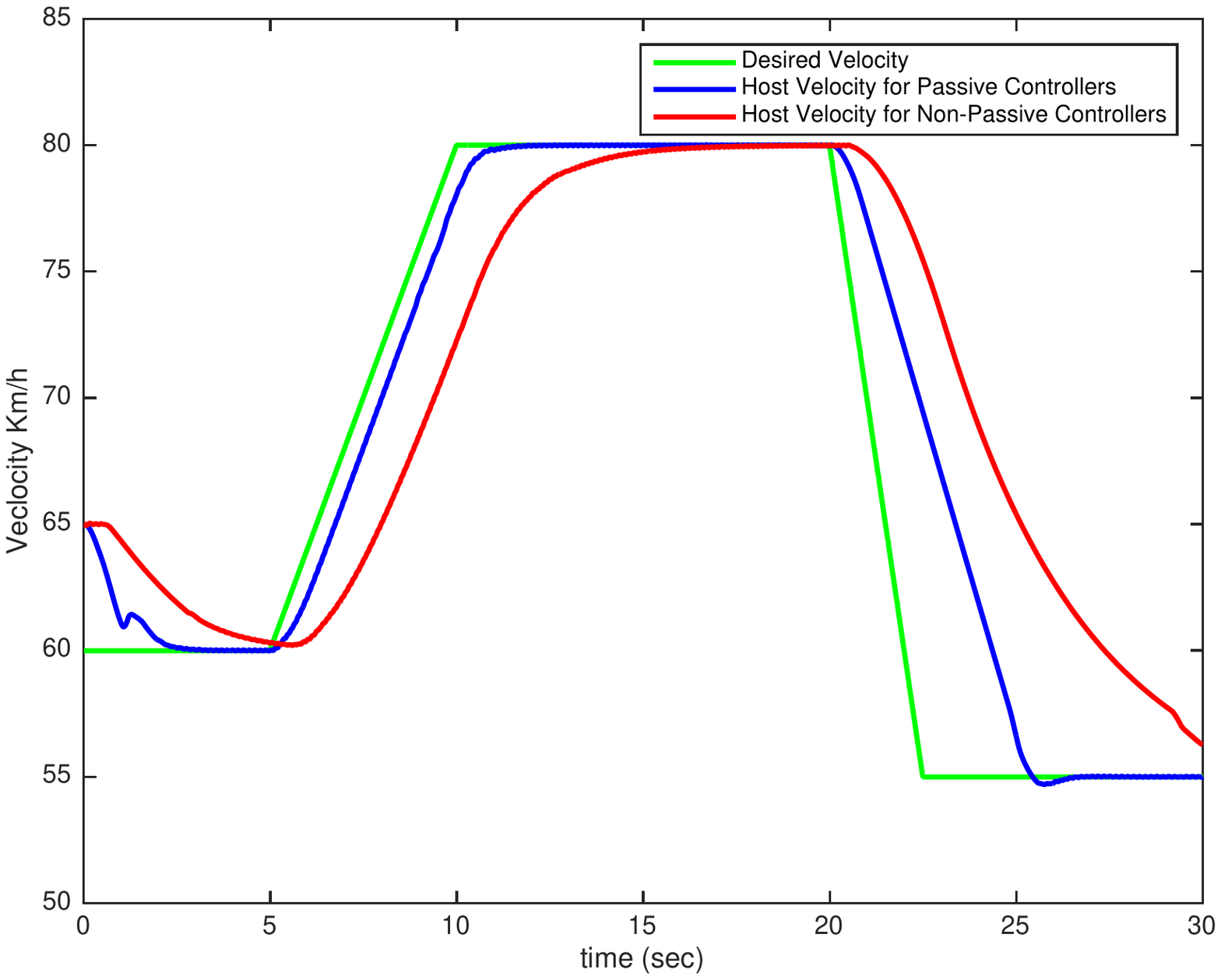}
		\caption{Performance Results for passive and nonpassive ACC for thirty seconds.}
		\label{fig:comp}
		\end{figure}

			Figure \ref{fig:sin} shows the simulation results for extremum-seeking control where the desired speed profile follows a sinusoidal structure. The design parameters used in this experiment are similar to the values given in Table 2.  Figure \ref{fig:Nyqsecc} shows the Nyquist and inverse Nyquist plots for this experiment that lie on the right-half plane and prove passivity for the controller. The difference in performance for passive and non-passive controllers becomes clearer in shorter experiments. Figure \ref{fig:comp} compares the performance for passive and non-passive controllers in the case where the experiment lasts for thirty seconds. 
			
				In our next experiment, we locate the lead vehicle close to the host vehicle, such that both the spacing control, and speed control will be activated at different points on the road depending on the relative distance between two cars. The desired velocity to be maintained by host vehicle with  initial speed of 60 km/h is the constant 80 km/h. The experiment runs for 120 seconds and the safe distance between two vehicles to be maintained is 10m. After the first 10s of our experiment, the radars on host vehicle detect the slower lead vehicle, and transit the controller into spacing control mode to control the safe distance between the cars. After 40s, the lead vehicle speeds up to 85 km/h and as a result, the velocity of the host vehicle increases to maintain the desired velocity, this scenario can be seen between the time segment 40-70s. For the time segment 70-90s, the lead vehicle slows down, and as a result the host vehicle also starts to decrease its velocity in order to maintain the safe distance between two vehicles. And at approximately 105s, two vehicles start to travel at the same speed again. In our following experiments we use the same speed profile (a user-defined trajectory) for the lead vehicle but for different time intervals $T$, so it would be beneficial to describe the lead vehicle velocity profile $v_l(t)$ in terms of $T$:
				
				\begin{eqnarray*}
					v_l(t) = \left\{ \begin{array}{lllll}
						60, \quad& \mbox{if ~~$ t\leq~\frac{1}{3}T$} \vspace{0.03in}\\
						60 + \frac{150}{T}(t-\frac{1}{3}T), \quad& \mbox{if~~ $ \frac{1}{3}T\leq~t \leq~\frac{1}{2}T$}\vspace{0.03in}\\
						85, \quad& \mbox{if ~~$ \frac{1}{2}T\leq~t \leq~\frac{2}{3}T$}\vspace{0.03in}\\
						55 + \frac{150}{T}(t-\frac{5}{6}T), \quad & \mbox{if~~ $ \frac{2}{3}T\leq~t \leq~\frac{5}{6}T$}\vspace{0.03in}\\
						55, \quad & \mbox{if ~~$ \frac{5}{6}T\leq~t \leq~T$}%\vspace{0.03in}\\
					\end{array} \right.
				\end{eqnarray*}

			Time-delays and actuator lag are the inevitable nature properties of the actuators and sensors in control systems. A random delay between 0.4 seconds to 0.6 seconds exists in each of PID controllers (spacing and velocity controllers) to account for existing delays in automotive control systems and also communication, computational, and human delays. Time-delays are inherent in many engineering systems, and they can degrade the performance of controlled systems, or cause instability. Stability conditions for PID controllers using Bifurcatino analysis has been obtained in \cite{silva2007pid}, but as the number of gain parameters increases, it becomes more difficult to use this approach. The act-and-wait control in which the inherently infinite structure of time-delay systems are described by finite-dimensional transition matrices is another possible solution for handling delays \cite{insperger2006act}. The problem of controlling a string of vehicles moving in one dimension and what is known as platoon stability of vehicles relies heavily on the performance of adaptive cruise controls implemented in vehicles in the presence of delays. Several approaches have been proposed in literature to overcome the effect of time delay in platoon stability \cite{guo2011hierarchical,yue2012guaranteed}. The effects of communication delay is investigated in \cite{guo2011hierarchical}. Additional to previous examples, linear models for human operators for drivers are also of the form input/output delays used in our simulation models. Our method provides another alternative to solve these problems. 
			
\begin{figure}[!t]
				\centering
				\includegraphics[scale = 0.65]{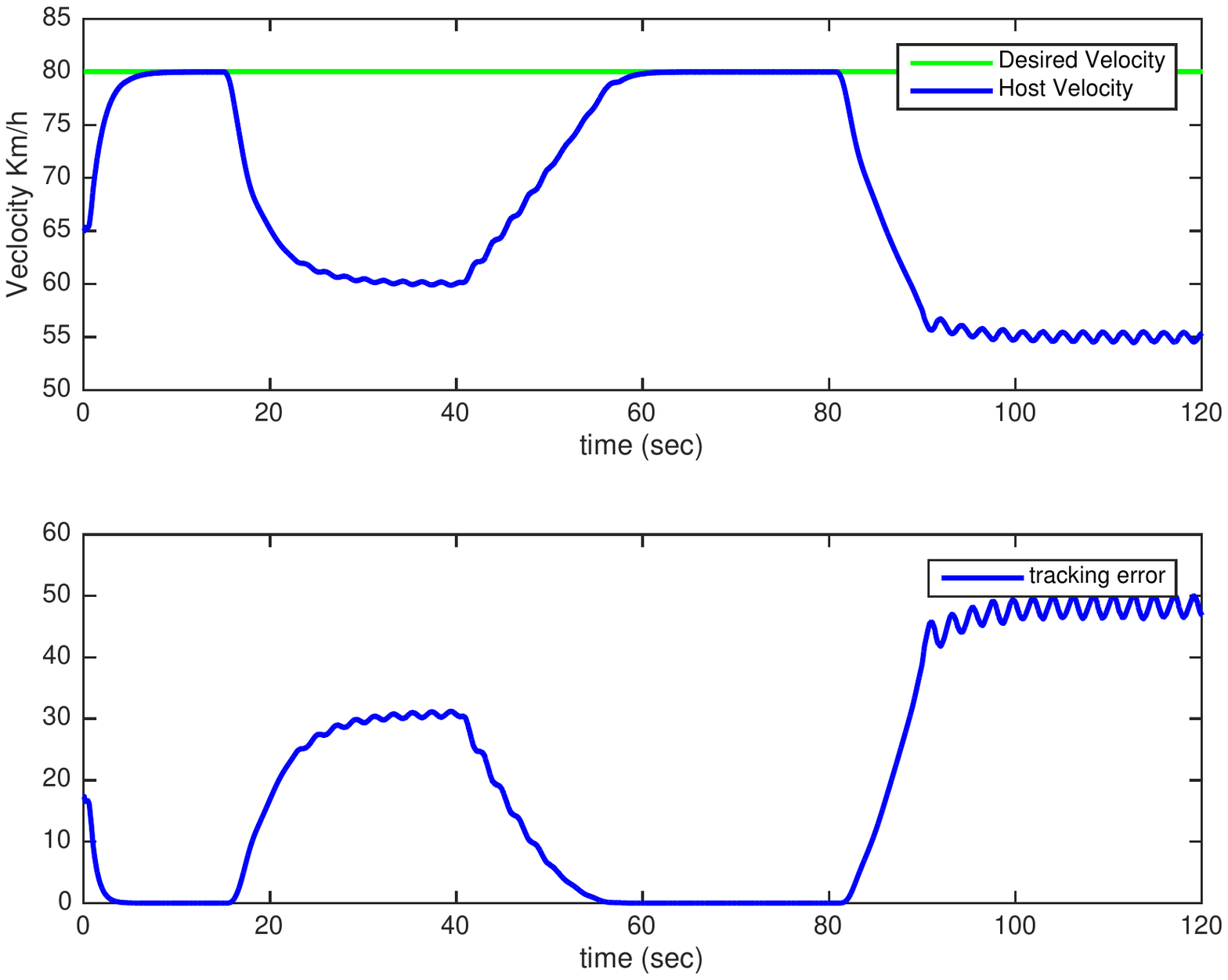}
				\caption{Simulation Results for Space/Speed Control With Delays.}
				\label{fig:ACCSpaceone}
				\end{figure}

			  Figure \ref{fig:ACCSpaceone} shows the results for the case where the system in the presence of delays is not passivated by our passivation methodology (the simulation results for the original non-passive set up with delays). As seen in the figures, the host velocity is not able to maintain the desired velocity for the entire course of experiment, given the fact that sometimes the gap between the host and lead vehicles is less than desired safe distance, hence the space control is turned on and the host velocity is forced to slow down. As shown in figure \ref{fig:ACCSpaceone}, the vehicle controlled by the non-passive controllers with delays suffers from a jittering performance caused by rapid breaking at different points of travel resulting in an uncomfortable ride experience.
				
				The adaptive cruise control in our last experiment includes a speed control side and a spacing control side, each containing delays. According to our method of passivation, two sets of $m$'s for each controller are used to passivate and overcome the negative effects of the delays and improve performance. Consequently, two different extremum-seeking schemes are designed to optimize these two sets of $m$'s. The design parameters given in table 2 are used for the extremum-seeking design that passivates the speed controller, and the design parameters given in table 4 are used to passivate the spacing controller.			
			
\begin{figure}[!t]
							\centering
							\includegraphics[scale = 0.65]{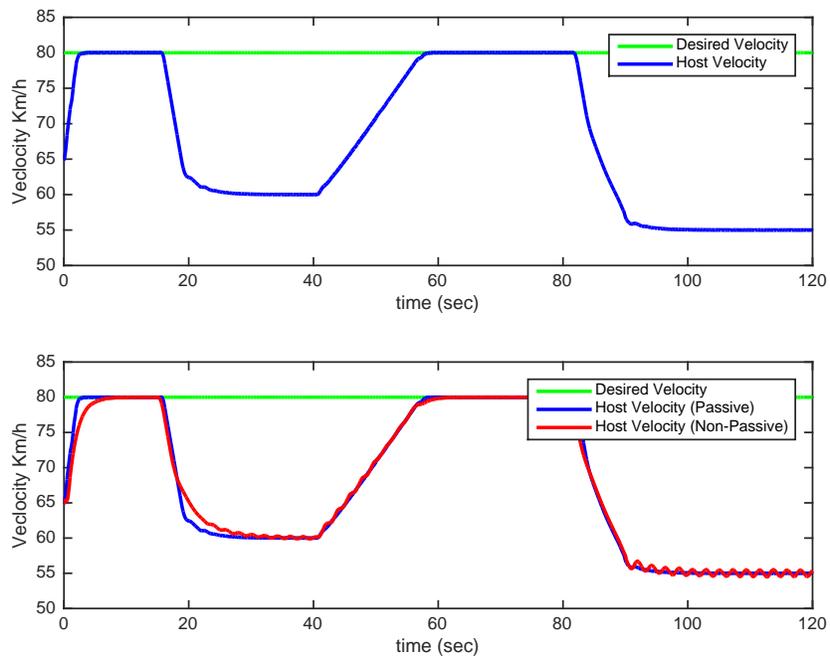}
							\caption{Simulation Results for Passive Space/Speed Control Set ups.}
							\label{fig:PassiveACC}
						\end{figure}

\begin{figure}[!t]
								\centering
								\includegraphics[scale = 0.65]{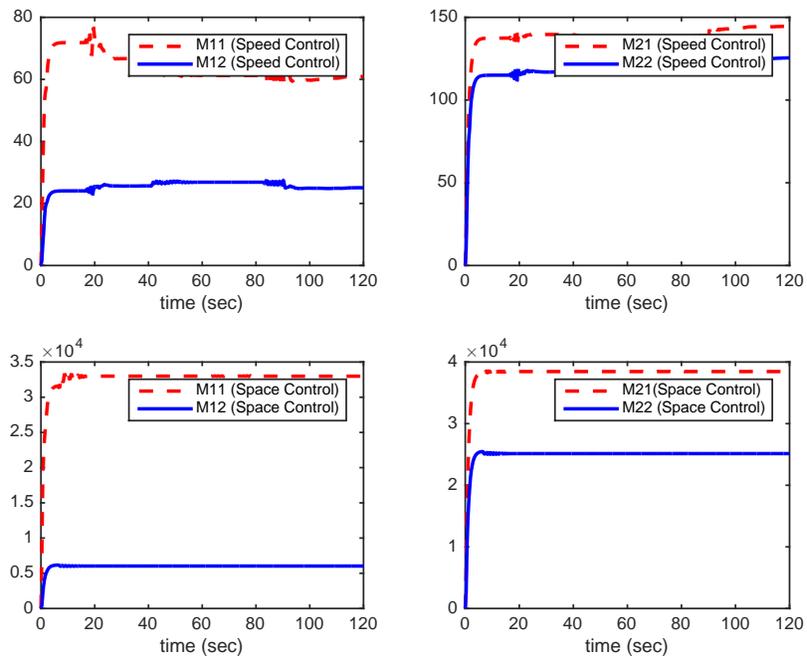}
								\caption{Simulation Results for Parameters for Passive Space/Speed Control Set up.}
								\label{fig:PassiveACCM}
							\end{figure}

\begin{figure}[!t]
		\centering
		\includegraphics[scale = 0.65]{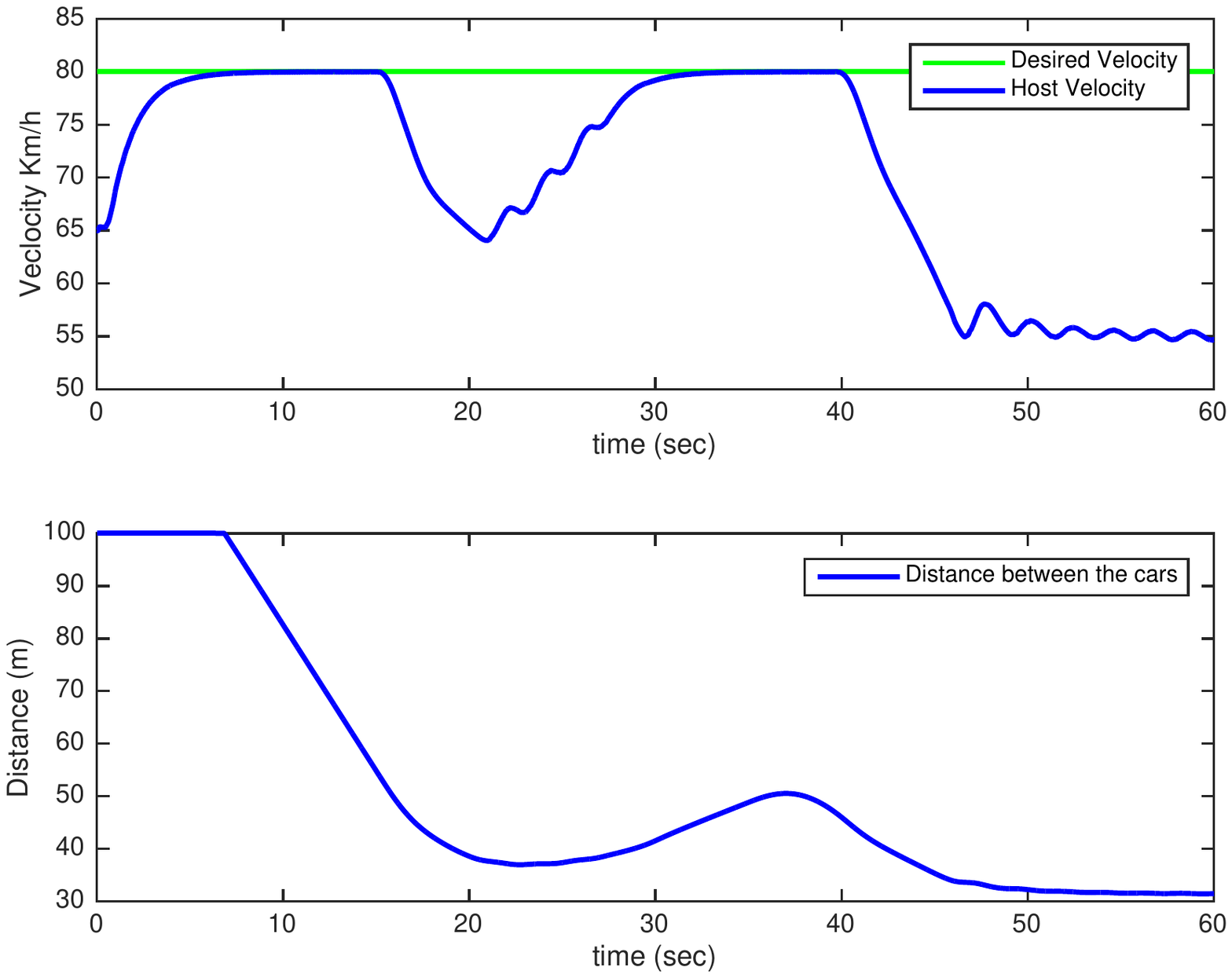}
		\caption{Simulation Results for Non-Passive Space/Speed Control Set up for 60s experiment.}
		\label{fig:NonPassiveshort}
	\end{figure}

\begin{figure}[!t]
									\centering
									\includegraphics[scale = 0.65]{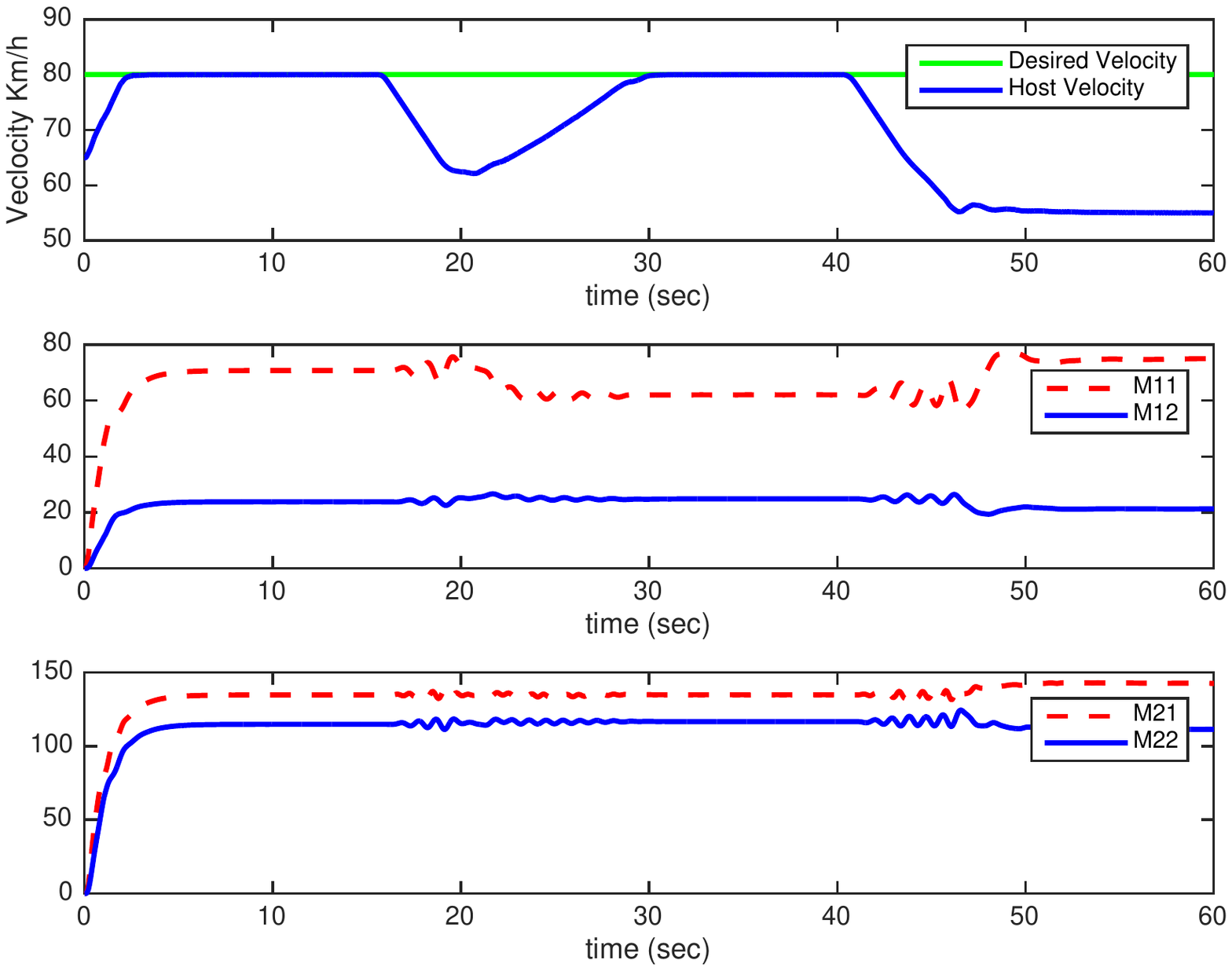}
									\caption{Simulation Results for Passive Space/Speed Control Set up for 60s experiment.}
									\label{fig:PassiveACCshort}
								\end{figure}

\begin{figure}[!t]
										\centering
										\includegraphics[scale = 0.65]{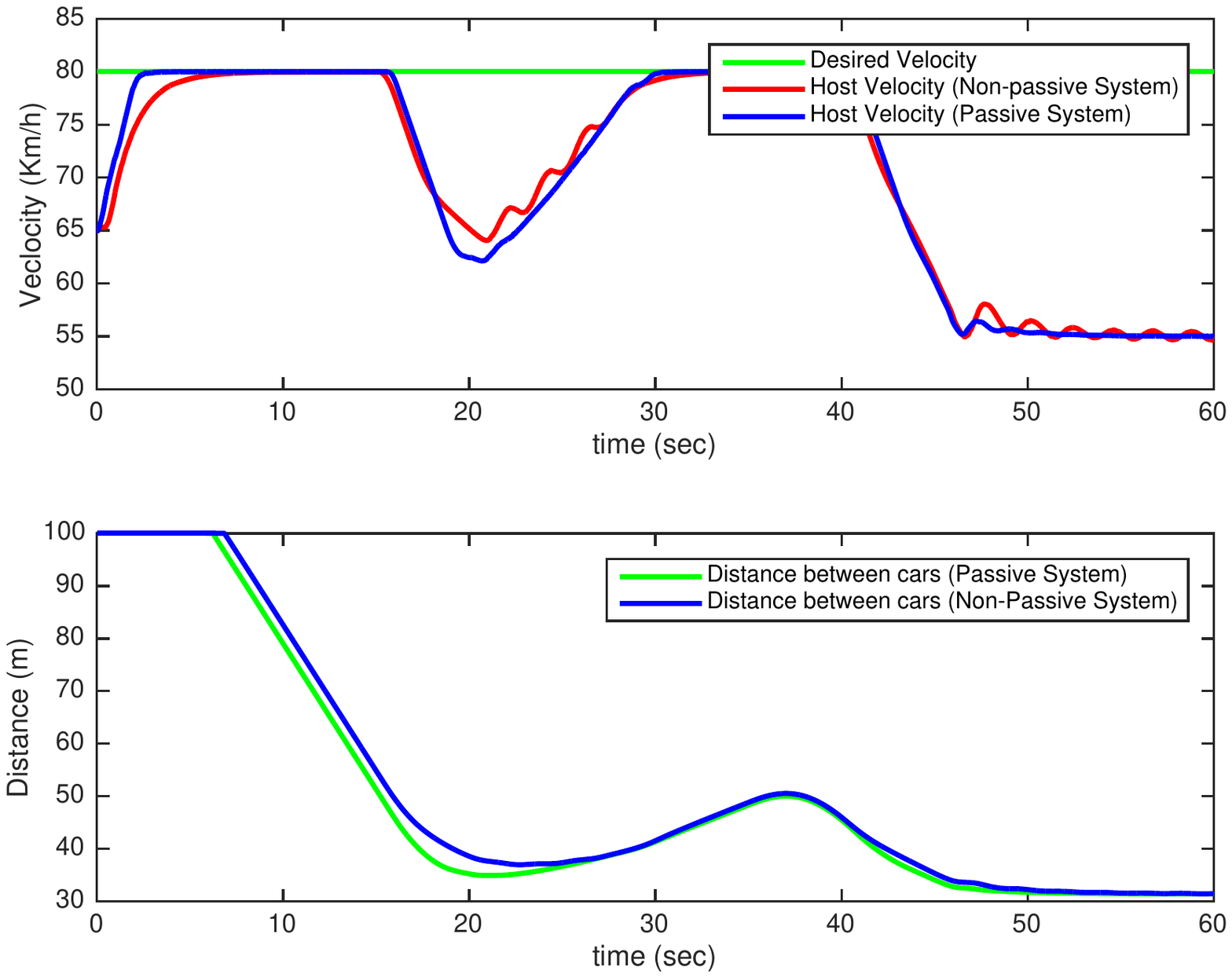}
										\caption{Passive Performance vs. Non-Passive Performance.}
										\label{fig:PassiveNonPassiveComp}
									\end{figure}

				Figure \ref{fig:PassiveACC} shows the simulation results for this experiment, and figure \ref{fig:PassiveACCM} depicts the values for the parameters that are being optimized. The simulation results prove that passivating the system enhances the overall performance, and overcomes the oscillations in host velocity resulting in less breaking, and a better ride experience for passengers. It is safe to say that the passive controllers give a more reliable overall presentation. The difference in performance between the passive and non-passive controllers becomes more dominant as we shorten the length of our experiment since the controller has less time to overcome the effects of the delays in shorter time intervals. Figure \ref{fig:NonPassiveshort}, and figure \ref{fig:PassiveACCshort} show the result for the same designs under a shorter experiment run-time of 60 seconds. As seen in figure \ref{fig:NonPassiveshort}, the non-passive controller gives a considerably less reliable performance, and the passive controller in figure \ref{fig:PassiveACCshort} is able to stabilize the system. Figure \ref{fig:PassiveNonPassiveComp} shows the performance for these two designs in one plot. Another interesting result to note is that the passive system makes use of the spacing control more efficiently, as seen in \ref{fig:PassiveACCshort}, the passive system gets closer to the lead vehicle while maintaining the safe distance, and the desired velocity for a longer time-interval. 
				
					\begin{table}[]
						\centering
						\caption{Design Parameters for Space Controller}
						\label{table4}
						\begin{tabular}{|c|c|c|}
							\hline
							& $\omega$  & $a$ \\  \hline 
							$m_{11}$  & $5+\frac{\pi}{2}$ & $0.01$ \\ \hline
							$m_{12}$  & $5$     & $0.01$\\ \hline
							$m_{21}$  & $10.14$       & $0.01$ \\ \hline
							$m_{22}$ & $7$        & $0.01$\\\hline
							
						\end{tabular}
					\end{table}
					
	%%===================================================================================								

			\section{Conclusion}
			
			Much effort has been invested in developing systematic tuning methods such as the ones given in surveys \cite{astrom1995theory,aastrom1993automatic}. But in an environment where the plant model is not known, extremum-seeking control as a real-time optimization technique combined with passivity is a great alternative. Compared to our previous experiments with the Hooke and Jeeves' optimization method, extremum-seeking control performs faster and requires only one experiment to achieve the optimal values - this decreases the time for finding correct passivation parameters from almost a quarter of an hour to seconds depending on the length of experiment. Some of other advantages of extremum-seeking control compared to Hooke and Jeeves' method are ESC's better robustness against disturbance, and noise, and its relatively less sensitive overall structure toward sudden changes. Additionally, extremum-seeking controller optimizes the system in real time. Our previous experiments also showed that one does not necessary need to know the shape of the tracking target in advance to select controller parameters accordingly, and that the optimized parameters obtained by our extreme-seeking co-simulation based framework platform for a certain experiment can still improve the behavior of system for different speed trajectories. The automotive control technology faces many control design challenges including system nonlinearities, insufficient modeling and measuring information, various sources of disturbances from outside, and constraints on efficiency and performance, and problems coming from scalability and interconnection between different units. In order to achieve a good controller system performance, the aforementioned control strategy meets all these challenges and can be seen as a promising technique for the future.

	%%===================================================================================

	\clearpage

	\bibliographystyle{IEEEtran}
	\bibliography{bibfile}
\end{document}